\documentclass[11pt]{amsart}

\usepackage{eucal,url,amssymb,stmaryrd,enumerate,amscd,verbatim}
\usepackage{amsmath}
\usepackage[pagebackref,colorlinks=true,linkcolor=blue,citecolor=blue]{hyperref}
\usepackage{amsfonts}
\usepackage{amsthm}
\usepackage[margin=1in]{geometry}
\usepackage{eurosym,lipsum,latexsym,esint,mathrsfs}
\numberwithin{equation}{section}
\newtheorem{thrm}{Theorem}[section]
\newtheorem{lemma}[thrm]{Lemma}
\newtheorem{prop}[thrm]{Proposition}
\newtheorem{cor}[thrm]{Corollary}

\newtheorem{rmrk}[thrm]{Remark}
\newtheorem{exm}[thrm]{Example}

\newcommand{\eps}{\varepsilon}

\newcommand{\R}{\mathbb{R}}
\newcommand{\Rn}{\mathbb{R}^n}

\newcommand{\Hn}{\mathbb{H}^n}

\newcommand{\bG}{\mathbb{G}}
\newcommand{\algg}{\mathfrak g}
\newcommand{\frlap}{\mathcal{L}_s}

\newcommand{\norm}[1]{\lVert#1\rVert}

\newcommand{\dom}{\mathcal{D}^{s,2}(\bG)}

\newcommand{\qn}[1]{\norm{#1}_{\,L^{2\text{*}(s)}(G)}}

\newcommand{\volg}{\, d{g}}
\newcommand{\volh}{\, d{h}}
\newcommand{\vol}{\, d{H}}
\newcommand{\bard}{\underline{d}}
\newcommand{\ddt}{\frac {d}{dt}}

\newcommand{\sA}{\mathscr A}

\def \de {\partial}

\begin{document}

\title[Fractional operators and Sobolev spaces on homogeneous groups]{Fractional operators and Sobolev spaces on homogeneous groups }

\begin{abstract}
    We establish foundational properties of fractional operators on Lie groups of homogeneous type. We prove embedding theorems for the associated Sobolev-type spaces.
\end{abstract}

\author{Nicola Garofalo}
\address[Nicola Garofalo]{School of Mathematical and Statistical Sciences\\ Wexler Hall, 901 Palm Walk Room 216\\ Tempe, AZ 85281, USA}
\vskip 0.2in
\email{nicola.garofalo@asu.edu}

\author{Annunziata Loiudice }

\address[Annunziata Loiudice]{Dipartimento di Matematica\\ Universit\`a degli Studi di Bari Aldo Moro\\Via Orabona, 4- 70125 Bari, Italy}

\email{annunziata.loiudice@uniba.it}

\author{Dimiter Vassilev}
\address[Dimiter Vassilev]{ Department of Mathematics and Statistics\\
University of New Mexico\\
Albuquerque, New Mexico, 87131-0001\\
}
\email{vassilev@unm.edu}

\date{\today}
\maketitle
\tableofcontents

\section{Introduction}
The general purpose of the present paper is to define and establish the basic properties of fractional operators and the corresponding Sobolev-type spaces on Lie groups of homogeneous type. By this denomination, we intend a connected and simply connected Lie group endowed with a non-isotropic family of dilations, in the sense of \cite{Stein68}, \cite{KS69}, \cite{KV71}, and \cite{FS82}. Our approach to introducing the relevant Sobolev-type space and its associated fractional operator is based entirely on an explicit Dirichlet form. For clarity of exposition, we restrict attention to the case which, in the Euclidean setting and up to a multiplicative constant, coincides with the fractional Laplacian $(-\Delta)^s$.

It is important to note that the Lie groups considered here are not assumed to be Carnot groups (also known as stratified nilpotent Lie groups). In particular, we do not impose the H\"ormander condition in \cite{Ho} on the invariant vector fields of homogeneity one. Nevertheless, even when these vector fields do not generate the entire Lie algebra, they play a basic role in the analysis. In this sense, the paper sheds new light on results concerning nonlocal operators in both Euclidean and Carnot group settings and provides a simple and unified proof of several foundational properties.
As an application of the present results, in the companion paper \cite{GLV22} we established Lebesgue-space regularity for an associated nonlocal Schr\"odinger-type equation, together with the sharp decay of nonnegative solutions corresponding to the extremals of the relevant Sobolev embedding theorem.
To state the main results, we now briefly recall the general setting, following \cite{FS82}.

\subsection{Notation and definitions}
Let $\bG$ be  a homogeneous group as defined in  \cite[Chapter 1]{FS82}. In particular, $\bG$ is a connected, simply connected nilpotent Lie group whose Lie algebra $\algg$ is endowed with a family of non-isotropic group dilations $\{\Delta_\lambda\}_{\lambda>0}$. This means that there is a basis $X_j$, $j=1,\dots,n$ of $\algg$, and positive real numbers $d_j$, such that,
\begin{equation}\label{e:dilations on algebra}
1=d_1\leq d_2\leq\dots\leq d_n \quad \text{and}\quad
\Delta_\lambda X_j=\lambda^{d_j} X_j.
\end{equation}
The exponential map $\exp:\algg\rightarrow \bG$ is a diffeomorphism onto. We can define a $1$-parameter family of automorphisms of the group $\bG$ by letting $\delta_\lambda = \exp\circ\Delta_\lambda\circ\exp^{-1}$. For convenience, we will henceforth use the same notation $\delta_\lambda$ for both the dilations on $\algg$ and $\bG$. As customary, we indicate by
\[
Q=d_1+\dots+d_n
\]
the homogeneous dimension of $\bG$.
Denoting with $\volg$ or $dH$ a fixed bi-invariant Haar measure on $\bG$, the following identity holds:  $d(\delta_\lambda g)=\lambda^Q\, dg$.

We shall assume that $\bG$ is endowed with a fixed homogeneous norm $|\cdot|$. This means that, if $e$ is the group identity, then $g\mapsto |g|$ is a function in $C(\bG)\cap C^\infty(\bG\setminus \{e\})$, such that $ |g|=0$ if and only if $g=e$. Furthermore, for every $g, h\in \bG$ one has
\begin{equation}\label{e:homog norm}
(i) \ |g^{-1}|=|g|; \quad (ii)\ |\delta_\lambda g|=\lambda|g|,
\end{equation}
and the triangle inequality
\begin{equation}\label{e:triangle ineqs}
 |g h|\leq |g|+|h|.
\end{equation}
We shall denote with $B_R(g)$ or $B(g,R)$ the open balls with center $g$ and radius $R$ defined using the fixed  homogeneous norm,
\begin{equation}\label{e:right-invarinat balls}
B_R(g)\equiv B(g,R)=\{h\, \mid\, |g^{-1}\, h|<R\}.
\end{equation}
For balls centered at the identity we will write  $B_R$, instead of $B_R(e)$.

We note that, according to \cite{HS90}, any homogeneous group supports a norm satisfying \eqref{e:triangle ineqs}. In Section \ref{ss:limits} we will also require certain symmetries of the fixed homogeneous norm. For example, the homogeneous norm constructed in \cite{HS90} is a poly-radial function of the coordinates \eqref{e:homog coord} grouped according to their homogeneity. However, in order to compute the limits of the fractional operator as $s\rightarrow0+$ and $s\rightarrow 1-$, we shall need a smaller group of  symmetries, which allows to use other homogeneous norms. In particular, in the case of the Heisenberg group, one can also use the Kor\'{a}nyi gauge.

If $L_{g_0}g=g_0\, g$ indicate the left-translations in $\bG$, the function $d(g,h)=|g^{-1} h|$ defines a distance on $\bG$, invariant with respect to $L_{g_0}$ for every $g_0\in \bG$. One has in fact
\begin{equation}\label{e:right inv dusatnce}
d(L_{g_0}g,L_{g_0}h)=|( g_0\, g)^{-1}\, ( g_0\, h)|=|g^{-1} g_0^{-1}\, g_0\,  h|=|g^{-1}\, h|=d(g,h),
\end{equation}
but the distance $d(g,h)$ is not necessarily invariant under right-translations. In particular, the balls \eqref{e:right-invarinat balls} are obtained by left-translating balls centered at the identity.
For $0<s<1$ let $\mathcal{D}^{s,2}(\bG)$ be the  fractional Dirichlet (homogeneous Sobolev) space defined as the closure of $C^\infty_0(\bG)$ with respect to the semi-norm
\begin{equation}\label{e:def frac norm}
[u]_{s,2}=\left(\int_{\bG}\int_{\bG}\frac {|u(g)-u(h)|^2}{| h^{-1}\, g|^{Q+2s}}\, {dg}dh\right)^{1/2}= \left(\int_{\bG}\int_{\bG}\frac {|u(g)-u(gh)|^2}{|h|^{Q+2s}}\, {dg}dh\right)^{1/2}.
\end{equation}
In particular, the $\mathcal{D}^{s,2}(\bG)$ norm is invariant under left-translations.

\subsection{Main results}
On the space $\dom$, we will use $\norm{u}_{\mathcal{D}^{s,2}(\bG)}$ for $
[u]_{s,2}$, hence the Dirichlet norm can be written as $$\norm{u}_{\mathcal{D}^{s,2}(\bG)}=\sqrt{\mathcal{D}_s(u,u)},$$
where $\mathcal{D}_s(u,\phi)$  is the quadratic form defined in \eqref{e: dirichlet form} below.

In Theorem \ref{t:trancation thrm} we sketch the usual identification of $\dom$ with  the space of functions $u\in L^{2^*(s)}(\bG)$ with finite $[\cdot ]_{s,2}$ semi-norm.

The infinitesimal generator of the quadratic form $\mathcal{D}_s(u,\phi)$ will be used to define the fractional operator $\frlap$. Our first result relates the definitions of a fractional operator defined by hyper-singular integral and the  Dirichlet form.

\begin{thrm}\label{p:symm of d form and hypersing int}

The quadratic form
\begin{equation}\label{e: dirichlet form}
\mathcal{D}(u,\phi)=\mathcal{D}_s(u,\phi)\overset{def}{=}\int_{\bG}\int_{\bG} \frac
{(u(g)-u(h))\left(\phi(g)-\phi(h)\right)}{| h^{-1}\, g|^{Q+2s}}\,
{dg}dh
\end{equation}
is a symmetric form on the space $C_0^\infty(\bG)$ of smooth functions with compact support. In fact, for $u,\, \phi\in C_0^\infty(\bG)$
we have
\begin{equation}\label{e:pv integral}
\mathcal{D}(u,\phi)=\int_{\bG}\phi(g)\, \frlap u (g) dg,
\end{equation}
where we have defined
\begin{equation}\label{e:frac L via 2nd order}
\frlap u(g)=\int_{\bG} \frac {2u(g)-u(gh)-u(gh^{-1})}{|h|^{Q+2s}}dh.
\end{equation}
Furthermore, we have
\begin{equation}\label{e:frac L via hypersing app app}
\frlap u(g)=2 \lim_{\eps\rightarrow
0}\int_{\bG\setminus B(g,\eps)} \frac {u(g)-u(h)}{|
h^{-1}\, g|^{Q+2s}}dh.
\end{equation}

\end{thrm}

Let us observe, cf. subsection \ref{ss:invariances frac operator}, that $\frlap$ commutes with left-translations, but not necessarily with right-translations.

When $\bG = \Rn$, up to a multiplicative constant, definition \eqref{e:frac L via 2nd order} gives the fractional Laplacian $(-\Delta)^s$, first introduced in \cite{Ri}, see also \cite{Stein61}. We recall in this connection that a useful alternative representation of this operator is obtained using the heat semigroup in $\Rn$, see e.g. \cite{G19}. Such representation continues to be valid whenever a good heat semigroup is available. For instance, in any Carnot group $\bG$ with sub-Laplacian $\mathcal L$, one can define
\begin{equation}\label{bala}
{\mathcal{L}}^s u(g) \overset{def}{=} (-\mathcal L)^s = - \frac{s}{\Gamma(1-s)} \int_0^\infty \frac{1}{t^{1+s}} \left(P_t u(g) - u(g)\right) dt,
\end{equation}
where $P_t = e^{-t\mathcal L}$, see \cite{Fo75}.
  However, the reader should be aware that in the non-Euclidean framework considered in the present paper, formula \eqref{bala} defines an operator that is markedly different from \eqref{e:frac L via 2nd order}. As this is one of the motivations for considering the general setting of homogeneous groups, we give some details in the case of groups of Heisenberg type. In fact, when $\bG$ is such a group, and the function
\[
|g| = N(z,\sigma) = (|z|^4 + 16 |\sigma|^2)^{1/4},\ \ \ \ g = (z,\sigma)\in \bG,
\]
denotes the Kor\'anyi-Folland gauge (see \cite{Ka} for the relevant definitions), then - up to a multiplicative constant -
$\frlap$ coincides with the \emph{geometric} fractional sub-Laplacian introduced in \cite{BFM}; see also \cite{FGMT} and \cite{RT}. The operator $\frlap$  was first introduced in \cite{BFM} via a spectral formula. This was generalized to groups of Heisenberg type $\bG$ in \cite{RT}. An important fact, first proved for $\Hn$ in \cite{FGMT} using hyperbolic scattering, and subsequently generalised to any group of Heisenberg type in \cite{RT} using non-commutative harmonic analysis, is the Dirichlet-to-Neumann characterisation of $\mathscr L_s$ using the solution to a certain extension problem from conformal CR geometry  different from that of Caffarelli-Silvestre in \cite{CS}. Yet another fundamental fact, proved in
\cite[Proposition 4.1]{RTaim} and \cite[Theorem 1.2]{RT} for $0<s<1/2$, is the following remarkable M. Riesz type representation
\begin{equation}\label{riesz}
 \alpha(m,k,s)\  \mathscr L_s u(g)=\int_{\bG} \frac {u(g)-u(h)}{|h|^{Q+2s}}dh.
\end{equation} Using the heat equation approach in \cite{GT20Adv}, \cite{GTinter}, see also \cite{GT23PA}, formula \eqref{riesz} was extended to  the whole range $0<s<1$. In \eqref{riesz} the number $\alpha(m,k,s)>0$ denotes an explicit constant depending on $s$ and the dimensions $m$ and $k$ of the horizontal and vertical layers of the Lie algebra of the group of Heisenberg type.
The conformal invariance of this pseudo-differential operator is manifested in the result proved in \cite{GT20Adv}, which establishes that its fundamental solution is given by
\[
\mathscr E_{(s)}(g) = C_{(s)}(m,k)\ N(z,\sigma)^{2s-Q},
\]
where $C_{(s)}(m,k)>0$ is an explicit universal constant. In contrast,
the fundamental solution of the non-geometric operator $\mathcal L^s$ in \eqref{bala}, also found in \cite{GT20Adv}, is given by the following expression
\begin{align*}
\mathscr E^{(s)}(g) & = C^{(s)}(m,k) \frac{1}{ |z|^{2(\frac m2 + k - s)}} \int_0^1 (\tanh^{-1} \sqrt y)^{s-1} \left(1-y\right)^{\frac m4-1}  y^{\frac 12(k - s-1)}
\\
& \times  F\left(\frac{1}{2}(\frac m2 + k - s),\frac 12(\frac m2 + k+1-s);\frac k2;- \frac{16|\sigma|^2}{|z|^4} y\right) \ dy,
\end{align*}
where $C^{(s)}(m,k)>0$ is another explicit universal constant, and $F(a,b;c;x)$ denotes Gauss' hypergeometric function. In particular, $\mathscr E^{(s)}(g)$ depends only on $|z|^4$ and $|\sigma|^2$, but it is not a function of the gauge.

Turning back to the general setting, the constants in front of the integrals defining the fractional operator are chosen so that formula \eqref{e:pv integral} holds true. A weak solution of the equation $\frlap u(g)=F$ is consistent with the usual formulation \[
\mathcal{D}(u,\phi)=\mathcal{D}_s(u,\phi){=}\int_{\bG}\phi (g)  F(g)\, dg.
\]
Properties of weak solutions of equations involving the fractional operator in the  general setting of homogeneous group were studied in \cite{GLV22}.

A consequence of Theorem \ref{p:symm of d form and hypersing int} and the definition of the Dirichlet space is that  the fractional operator $\frlap$ defines a self-adjoint operator on $L^2(\bG) $, which can be thought of as the infinitesimal generator of the associated diffusion.

The following result, whose proof can be found in Section \ref{ss:limits}, can be considered as a version of the Euclidean \cite[Proposition 4.4]{dNPV12}. We denote by $d\sigma$ the surface measure appearing in the integration in polar coordinates formula \eqref{e:polar coordinates integral}.

\begin{thrm}\label{p:limits frac operator}
Let  $m$ be the number of left-invariant vector fields $X_j$ of homogeneity one with associated coordinates $x_j$, cf. \eqref{e:homog coord}. Suppose the fixed  homogeneous norm $|g|$ is a radial function of the coordinates $(x_1,\dots,x_m)$, i.e., its dependence on these coordinates is only as a smooth function of $x_1^2+\dots+x_m^2$ on $x_1^2+\dots+x_m^2\not=0$.
For a function $u\in C^\infty_0(\bG)$ we have the identities
\begin{equation}\label{e:limits frac operator}
\lim_{s\rightarrow 0+} \frac {s}{\sigma_Q}\frlap u(g)=u(g) \quad\text{and}\quad \lim_{s\rightarrow 1-} \frac {2m(1-s)}{\tau_m}\frlap u(g)=-\sum_{i=1}^m X_i^2u(g),
\end{equation}
where $\sigma_Q$ is the measure of the unit sphere (defined by the fixed norm) and
\begin{equation}\label{e:tau_m def}
\tau_m=\int_{|h|=1} {\sum_{i=1}^m |x_i(h)|^2}\, d\sigma(h).
\end{equation}
\end{thrm}

With the help of Theorem \ref{p:limits frac operator} we are able to establish the following connection between the fractional Dirichlet norm and the $L^2$ Dirichlet norm of a smooth function with compact support inspired by the Euclidean results of Bourgain, Brezis and Mironescu \cite{BBM01}, and Maz'ya and Shaposhnikova \cite{MS02}. In the framework of Carnot groups, universal characterizations of BV and Sobolev spaces (with dimension independent constants) were recently obtained using heat semigroup techniques in \cite{GT24, GT23CAG}, see also \cite{BGT22} and \eqref{ms} below, for a related result. In the same framework, we mention that results such as Theorem \ref{p:limits frac operator} were obtained in \cite[Theor. 1.1, Cor. 1.1 \& Sec. 4]{CMSV21}, where the authors consider Orlicz-Sobolev spaces and \cite[Theor. 3.6]{Barb11}.  Using interpolations techniques, asymptotic versions of \eqref{1} and \eqref{2} below,  in the Euclidean and Carnot group settings, respectively, were obtained in in \cite{Mil05} and  \cite{MP19}.
Further developments will be discussed after the proof of Proposition \ref{t:BBM and MS} where, with the assumptions of Theorem \ref{p:limits frac operator}, we show the following results:
\begin{equation}\label{1}
\lim_{s\rightarrow 1-} \frac {2m(1-s)}{\tau_m}\norm{u}^2_{\mathcal{D}^{s,2}(\bG)}=\int_{\bG}\sum_{i=1}^m |X_iu(g)|^2 dg.
\end{equation}
and
\begin{equation}\label{2}
\lim_{s\rightarrow 0+} \frac {s}{\sigma_Q}\norm{u}^2_{\mathcal{D}^{s,2}(\bG)}=\int_{\bG} |u(g)|^2 dg.
\end{equation}
It is worth stressing that the sum of squares of vector fields that appears in the right-hand side of \eqref{e:limits frac operator} is, in general, not a H\"ormander type operator as we do not assume that the vector fields of homogeneity one satisfy a finite step generating condition. A non-Carnot group  example, which demonstrates this fact and illuminates further the discussions above is given in Example \ref{ex:non carnot}.

In the final section of the paper, we give the basic embedding results for the Dirichlet space $\mathcal{D}^{s,2}(\bG) $ associated to the quadratic form $\mathcal{D}$ given in \eqref{e: dirichlet form}, also called the homogeneous fractional Sobolev space of order $s$. This is the space  defined as the completion of $C_0^\infty(\bG)$ with respect to the norm \eqref{e:def frac norm}. The space $\mathcal{D}^{s,2}(\bG) $ is a space of functions in $L^{2^*(s)}(\bG)$ due to the  fractional Sobolev inequality in the following theorem. Furthermore, the norm of the embedding is achieved. The extremal functions satisfy a Yamabe type equation, which in the elliptic and sub-elliptic settings have drawn considerable attention. The characterization of all non-negative solutions of the obtained fractional Yamabe type equation is in general open,  a particular case being the fractional Yamabe equation on the Heisenberg (type) groups considered in  \cite{FGMT, RT, GTinter}.
We also give a Rellich-Kondrachov compact embedding result encompassing all of the known results in ''flat'' Euclidean or Carnot group settings.

\begin{thrm}\label{t:frac homog Sobolev and Rellich embed}
Let $Q$ denotes the homogeneous dimension of the  group $\bG$, $0<s<1$.

a) For some constant $S=S(s,Q)$ we have  for any $u\in C^\infty_0(\bG)$ the inequality
\begin{equation}\label{e:fracsobolev}
\left ( \int_{\bG} |u|^{2^*(s)}(g){dg}\right)^{1/2^*(s)}\leq S
\left(\int_{\bG}\int_{\bG} \frac {|u(g)-u(h)|^2}{|
h^{-1}\, g|^{Q+2s}}\, {dg}dh \right )^{1/2}
\end{equation}
where   the $L^2$ fractional Sobolev exponent is
\begin{equation}\label{e:r def}
2^*(s)\overset{def}{=}\frac{2Q}{{Q-2s}},
\end{equation}

b) The equality in the fractional Sobolev inequality \eqref{e:fracsobolev} on $\bG$ is achieved.

c) If $\Omega$ is a bounded domain, then the inclusion $\iota:\mathcal{D}^{s,2}(\bG) \hookrightarrow L^q(\Omega)$ is compact for any $q\in [1,2^*(s))$.
\end{thrm}

The proofs of the above results are given in Section \ref{ss:frac homog Sobolev and Rellich embed}, Theorem \ref{t:existence of extremals},  and Theorem \ref{t:frac Rellich embed}.
Let us note that for $h\in \bG$ and $\lambda >0$ letting
\begin{equation}\label{e:scaling}
L_hu=u\circ L_h \qquad \text{and}\qquad u_\lambda\ \equiv\ \lambda^{-Q/2^*(s)}\  \delta_{1/\lambda} u\ \overset{def}{=}\ \lambda^{-Q/2^*(s)}\ u\circ \delta_{1/\lambda},
\end{equation}
a small calculation shows that the norms involved in the fractional Sobolev inequality are invariant under the above defined left-translations and rescalings \eqref{e:scaling}. The $L^{2^{*}(s)}(\bG)$ norm is invariant also with respect to right-translations, but the $\mathcal{D}^{s,2}(\bG)$ norm is invariant only by the left-translations.

\subsubsection*{\textbf{Acknowledgements}} D.V. would like to thank the Department of Mathematics of Universit\`a degli studi di Bari Aldo Moro for the hospitality and resources during his visits. Funding for D.V. was partially provided by the the Efroymson Fund of the Department of Mathematics and Statistics, University of New Mexico. A.L. is a member of the \emph{Gruppo Nazionale per
l'Analisi Mate\-matica, la Proba\-bilit\`{a} e le loro
Applicazioni} (GNAMPA) of the \emph{Istituto Nazionale di Alta
Matema\-tica} (INdAM) and she is supported by the INdAM - GNAMPA
Project 2025 - CUP E5324001950001: \textit{Aspetti qualitativi per
equazioni nonlineari subellittiche su gruppi di Lie}.

\section{Background}

This sections serves to collect some material which is used in the rest of the paper.
\subsection{The Haar measure} Without loss of generality we will assume that the fixed Haar measure is the push forward of the Lebesgue measure on the Lie algebra via the inverse of the exponential map, see \cite[Proposition 1.2]{FS82}. We note that this gives a bi-invariant Haar measure.
The polar coordinates formula for such measure gives the existence of a unique Radon measure $d\sigma(h)$, such that, for $u\in L^1(\bG)$  we have the  identity,  \cite[Proposition (1.15)]{FS82} and \cite[p. 604]{KV71},
\begin{equation}\label{e:polar coordinates integral}
\int_{\bG} u(h)dh=\int_0^\infty\int_{|h|=1} u(\delta_r h) r^{Q-1}d\sigma(h)dr.
\end{equation}
In particular, we have, see \cite{FR66},
\begin{equation}\label{e:polar coord}
\int_{r<|g|<R}|g|^{-\gamma}\volg=\left\{
                  \begin{array}{ll}
                    \frac {\sigma_Q}{Q-\gamma}\left(R^{Q-\gamma} - r^{Q-\gamma} \right), & \gamma\not=Q\\
                    &\\
                    \sigma_Q\log(R/r), & \gamma=Q
                  \end{array}
                \right.
\end{equation}
for  a certain constant $\sigma_Q$ equal to the measure of the unit sphere (defined by the fixed norm),
\[
\sigma_Q=\int_{|h|=1} d\sigma(h).
\]

\subsection{Taylor formulas}

Here, we recall
the left Taylor formula and inequalities for  $u(gh)-u(g)$ attached to the considered homogeneous group relying on  \cite[Lemma 5.2]{KV71},  \cite[Theorem (1.37)]{FS82} and \cite[Theorem 2]{Bon09}.

Let $X_1,\dots X_n$ be a basis of the Lie algebra $\algg$ consisting of eigenvectors of the dilations $\delta_\lambda$ with corresponding eigenvalues $\lambda^{d_1}, \dots, \lambda^{d_n}$. We will denote with $X_j$, resp. $Y_j$, $j=1,\dots, n$ the corresponding left-invariant, resp. right-invariant, vector fields defined as usual by their action on a smooth function $f$ on $\bG$ by the formulas
\begin{equation}\label{e:inv vec fields}
(X_jf)(g)=\ddt f(g\exp(tX_j))\vert_{t=0} \quad \text{and}\quad (Y_jf)(g)=\ddt f(\exp(tX_j)g)\vert_{t=0}.
\end{equation}
The translation invariance implies that each of the vector fields $X_j$ and  $Y_j$ is formally skew-symmetric \cite[p. 21]{FS82}, i.e., for $u,\, v\in C^\infty_0(\bG)$ and $j=1,\dots,n$ we have
\begin{equation}\label{e:skew-adjoint v.f.}
\int_{\bG} u(g)X_jv(g)\, dg=-\int_{\bG} v(g)X_j u (g)\, dg\quad\text{and}\quad \int_{\bG} u(g)Y_jv(g)\, dg=-\int_{\bG} v(g)Y_j u (g)\, dg.
\end{equation}
The canonical coordinate functions of the first kind on $\algg$ with respect to the fixed basis $X_1,\dots X_n$ will be denoted by $x_1,\dots,x_n$, thus for $x\in\algg$ the equation $x=(x_1,\dots,x_n)$ means
\begin{equation}\label{e:homog coord}
x=x_1X_1+\dots + x_nX_n.
\end{equation}
We define the Euclidean norm  $\norm{\, .\,}$ on $\algg$ by declaring the fixed basis $X_1,\dots X_n$ to be orthonormal,
$$\norm{x}=\left(\sum_{j=1}^n x_j^2\right)^{1/2}.$$ In particular, using the exponential map, we obtain also a function on $\bG$, by defining $\norm{g}=\norm{\exp^{-1} g}$. In the coordinates $(x_1,\dots,x_n)$ on $\algg$ we have $X_j\vert_e=Y_j\vert_e=\partial_{x_j}\vert_0$ and the dilations take the form, recall the normalization $d_1=1$,
\[
\delta_\lambda(x)=(\lambda x_1,\lambda^{d_2}x_2\dots,\lambda^{d_n}x_n).
\]
We will refer to these coordinates also as the homogeneous coordinates.

For $k\in \mathbb{N}$ and a multi-index $I=(i_1,i_2,\dots,i_k)$ with $i_1, i_2, \dots,i_k\in \{1,2,\dots,n\}$ we let
\[
X_I=X_{i_1}\cdot\cdot\cdot X_{i_k} \quad\text{and}\quad x_I(g)=x_{i_1}(g)\cdots x_{i_k}(g).
\]
We denote with $|I|=k$ the order of the differential operator $X_I$, while
\[
d(I)=d_{i_1}+\dots+d_{i_k}
\]
will stand for the homogeneous degree or the degree of homogeneity of the operator with respect to the fixed dilations. We shall denote with $\Delta$ the additive semi-group generated by $0,d_1,\dots,d_n$. For $a\in \Delta$, due to  the condition $d_1=1$, the integer part $ [a]$  of $a$ satisfies trivially the formula, cf. \cite[Lemma (1.36)]{FS82},
\[
[a]=\max \{ |I|\, :\, d(I)\leq a \}.
\]

According to  \cite[Theorem (1.37)]{FS82}, see also \cite[Theorem 2]{Bon09},  there exists a constant $c$ depending on $\bG$ and the fixed norm, such that, for $a\in \Delta$ the Taylor expansion of order $a$ centered at $g$ of a smooth function near $g$ is given by
\begin{equation}\label{e:general Taylor}
u(gh)=P_{a}(h;g)\ +\ R_{a}(h;g),
\end{equation}
where the left Taylor polynomial of $u$ centered at $g$ of degree $a$ is defined by the formula
\begin{equation}\label{e:Taylor pol ord bard}
P_{a}(h;g)=u(g)+\sum_{I:\,  |I|\leq [a], d(I)\leq a} \frac {1}{|I|}(X_I u)(g)\,x_I(h),
\end{equation}
and the remainder $R_{a}(h;g)$ satisfies the estimate
\begin{multline}\label{e:Taylor rem ord a}
|R_{a}(h;g)|\leq \sum_{k=1}^{[a]+1}\frac {c^k}{k!}\sum_{I:\, |I|=k, d(I)> a} |h|^{d(I)} \, \sup_{|h'|\leq c|h|} |X_Iu(gh')|\\
 \leq  e^c\sum_{k=1}^{[a]+1}\sum_{I:\, |I|=k, d(I)> a} |h|^{d(I)} \, \sup_{|h'|\leq c|h|} |X_Iu(gh')| .
\end{multline}

\begin{rmrk}\label{r:left-inv Taylor}
We note that the Taylor polynomials for the left translation $u(hg)-u(g)$ use the right-invariant vector fields, cf. \cite[Theorem (1.37) and Remark p.31]{FS82}.
\end{rmrk}

It will be useful to write explicitly the Taylor formula of order two that we will use.
For $a=2$, noting that $2\in \Delta$ due to the assumption $d_1=1$, we obtain
\begin{equation}\label{e:Taylor pol ord 2}
P_{2}(h;g)=u(g)+\sum_{i:\, d_i\leq 2} (X_iu)(g)\,x_i(h)\ +\  \frac 12\sum_{i,j:\, d_i=d_j=1} (X_iX_ju)(g)\,x_i(h)x_j(h).
\end{equation}
The remainder $R_{2}(h;g)$ satisfies the estimate
\begin{equation}\label{e:Taylor rem ord 2}
|R_{2}(h;g)|\leq \sum_{k=1}^{3}\frac {c^k}{k!}\sum_{I:\, |I|=k, d(I)> 2} |h|^{d(I)} \, \sup_{|h'|\leq c|h|} |X_Iu(gh')|.
\end{equation}
It follows that there exist constants $d>2$, equal to the smallest $d(I)>2$ appearing in the above identity,  and $T_2=e^c>0$ such that for $|h|\leq 1$ we have
\begin{equation}\label{e:Taylor rem ord 2 bis}
|R_{2}(h;g)|\leq T_2 |h|^{d} \sum_{I:\, 1\leq |I|\leq 3, d(I)> 2} \sup_{|h'|\leq c|h|} |X_Iu(gh')|.
\end{equation}
Furthermore, we have
\begin{equation}\label{e:2nd order difference}
u(gh)+u(gh^{-1})-2u(g)=  \sum_{i,j:\, d_i=d_j=1} (X_iX_ju)(g)\,x_i(h)x_j(h) + \mathcal{R}(h;g),
\end{equation}
where from \eqref{e:Taylor rem ord 2 bis} we have for some $d>2$ the estimate
\begin{equation}\label{e:2nd order difference remainder}
\vert \mathcal{R}(h;g) \vert \leq 2T_2 |h|^{d} \sum_{I:\, 1\leq |I|\leq 3, d(I)> 2} \sup_{|h'|\leq c|h|} |X_Iu(gh')|, \qquad |h|\leq 1.
\end{equation}
 Indeed, on one hand  the homogeneous norm satisfies $|h^{-1}|=|h| $. On the other hand,  by the Campbell-Baker-Hausdorff formula  $x_i(h^{-1})=-x_i(h)$. Thus, an application of  \eqref{e:Taylor pol ord 2} and\eqref{e:Taylor rem ord 2} gives the claimed identity \eqref{e:2nd order difference}.

We state explicitly another corollary of the Taylor expansion, namely the mean value theorem \cite[(1.33)]{FS82}, according to which there exist constants $T_0$ and $c$, such that, for $u\in \mathcal{C}^1(\bG)$ and $g,\, h\in\bG$ we have
\begin{equation}\label{e:MVT}
\vert u(gh)-u(g)\vert \leq T_0\sum_{j=1}^n |h|^{d_j}\sup_{|h'|\leq c|h|} |X_ju(g h')|,
\end{equation}
noting explicitly that the $\sup$ involves derivatives along \emph{all} vector fields. We also recall \cite[Lemma 5.2]{KV71} for another version for homogeneous functions in  a more general setting.

\section{The Dirichlet form and fractional operator}\label{s:frac operator}

The proof of Theorem \ref{p:symm of d form and hypersing int} will use several estimates of integrals that we state as separate lemmas below. It is worth comparing the proofs of the analogous results in the Euclidean case, see \cite{Silvestre07},  \cite{dNPV12} and \cite{G19}, with the subtleties of the arguments   in the considered setting of non-commutative vector fields and variables of varied homogeneity.

We begin by extending the fractional operator to a larger space, including all constant functions, on which the fractional operator will vanish, and smooth functions that are  in $L^p(\bG)$, $1\leq p\leq \infty$.
First, we introduce the space
\[
L_s^p(\bG)=\{ u\in L^1_{loc}(\bG) \mid \int_{\bG} \frac {|u(h)|^{p-1}}{1+|h|^{Q+2s}} \, dh  <\infty\},
\]
see \cite{La}.
We will write $L_s(\bG)$ for the space $L_s^2(\bG)$.

\begin{prop}\label{p:defs of frac operator}
Let $u\in L_s(\bG)\cap C ^\infty(\bG)$. The integrals \eqref{e:frac L via 2nd order} and \eqref{e:frac L via hypersing app app} are finite and equal to each other.
\end{prop}

\begin{proof}
For the first part of the statement, it suffices to show that for a fixed $g\in\bG$ the following function of $h$ is integrable,
\begin{equation}\label{e:frac quotient in L1}
\frac {2u(g)-u(gh)-u(gh^{-1})}{|h|^{Q+2s}}\in L^1(\bG).
\end{equation}
Since the claim is trivial for $g=e$, we assume $|g|>0$.
Recalling \eqref{e:right-invarinat balls}, consider separately the integrals of \eqref{e:frac quotient in L1} on $B_1$ and its complement.
On $\bG\setminus B_1$ we rely on the condition  $u\in L_s(\bG)$.
 Using  the triangle inequality and the change of variables $h'=gh$ and $h'=h^{-1}$, we have
\begin{equation}\label{e:2nd order diff ext B1}
\int_{\bG\setminus B_1}\frac {\left| 2u(g)-u(gh)-u(gh^{-1}) \right|}{|h|^{Q+2s}}\, dh \leq   2\int_{\bG\setminus B_1} \frac {\left|u(g)\right|}{|h|^{Q+2s}}\, dh + 2\int_{\bG\setminus B_1(g)} \frac {\left|u(h)\right|}{|g^{-1}h|^{Q+2s}}\, dh.
\end{equation}
Notice that for $|h|\geq 1$ we have trivially $1+|h|^{Q+2s} \leq 2|h|^{Q+2s}$. The first integral in the right-hand side is then clearly finite by  $u\in L^1(\bG, \frac {dh}{1+|h|^{Q+2s}} )$. To see the finiteness of the second integral, let $2R=|g|$ and majorize as follows
\begin{multline}
\int_{\bG\setminus B_1(g)} \frac {\left|u(h)\right|}{|g^{-1}h|^{Q+2s}}\, dh \leq  \int_{(\bG\setminus B_1(g))\cap B_R} \frac {\left|u(h)\right|}{|g^{-1}h|^{Q+2s}}\, dhh + \int_{(\bG\setminus B_1(g))\setminus B_R} \frac {\left|u(h)\right|}{|g^{-1}h|^{Q+2s}}\, dh\\
\end{multline}
The first integral the right-hand side is finite since the domain is bounded and the integrand is a continuous function.  On the other hand, since $|g^{-1}h|\geq |h|-|g|\geq |h|/2$ on the complement of $B_R$, the last integral is finite by $u\in L_s(\bG)$.

On the other hand, for $|h|\leq 1$ we have by \eqref{e:2nd order difference} and  \eqref{e:2nd order difference remainder} the estimate
\begin{multline}\label{e:2nd difference quot}
 \left| \frac {2u(g)-u(gh)-u(gh^{-1})}{|h|^{Q+2s}}\right|
\leq  \sum_{i,j:\, d_i=d_j=1}\left| (X_iX_ju)(g)\right|\,  \frac{\left| x_i(h)x_j(h)\right|}{|h|^{Q+2s}}  \  + \   \frac {\left| \mathcal{R}(h;g) \right|}{|h|^{Q+2s}}\\
\leq  2C(T_2+1)  \sup_{B_{c}(g)} \{ \left |X_Iu \right| \  :  \   1\leq |I|\leq 3, \, 2\leq d(I)  \} \ \frac{|h|^{2}}{|h|^{Q+2s}},
\end{multline}
with a  constant $C=C(Q)$, noting that $|g^{-1}gh'|\leq |h'|\leq c|h|\leq c$ in \eqref{e:2nd order difference remainder},  and  that $d>2$ hence $|h|^d\leq |h|^2$ on $B_1$. An integration over the unit ball and  \eqref{e:polar coord} give
\begin{multline}
 \int_{B_1}\left| \frac {2u(g)-u(gh)-u(gh^{-1})}{|h|^{Q+2s}}\right| dh \\
\leq   2C(T_2+1)  \sup_{B_{c}(g)} \{ \left |X_Iu \right| \  :  \   1\leq |I|\leq 3, \, 2\leq d(I)  \}  \,\frac {\sigma_Q}{2(1 -s)}.
\end{multline}
The proof of \eqref{e:frac quotient in L1}, and therefore \eqref{e:frac L via 2nd order}, is complete.

We turn to the proof of \eqref{e:frac L via hypersing app app}, i.e., the equivalence of the two definitions of the fractional operator . Let $\eps>0$. From the translation invariance of the Haar measure, using the substitutions $h'= h^{-1}\, g$ and $h'^{-1}= h^{-1}\, g$ we have

\[
\int_{\bG\setminus B(g,\eps)} \frac {u(g)-u(h)}{|
h^{-1}\, g|^{Q+2s}}dh=\int_{\bG\setminus B_\eps} \frac {u(g)-u(gh^{-1})}{|h|^{Q+2s}}dh=\int_{\bG\setminus B_\eps} \frac {u(g)-u(gh)}{|h|^{Q+2s}}dh,
\]
after taking into account that $|h'^{-1}|=|h'|$. Therefore,
\begin{equation}\label{e:pv to non pv}
\int_{\bG\setminus B(g,\eps)} \frac {u(g)-u(h)}{|
h^{-1}\, g|^{Q+2s}}dh =\frac 12\int_{\bG\setminus B_\eps} \frac {2u(g)-u(gh)-u(gh^{-1})}{|h|^{Q+2s}}dh.
\end{equation}
From \eqref{e:frac quotient in L1} and Lebesgue's dominated convergence theorem we have
\[
\lim_{\eps\rightarrow 0} \int_{\bG\setminus B(g,\eps)} \frac {u(g)-u(h)}{|
h^{-1}\, g|^{Q+2s}}dh = \frac 12\int_{\bG} \frac {2u(g)-u(gh)-u(gh^{-1})}{|h|^{Q+2s}}dh,
\]
which completes the proof.
\end{proof}

We emphasize that, since $\frac {1}{1+|h|^{Q+2s}}\in L^1(\bG )$ by \eqref{e:polar coord}, it follows that $L^\infty(\bG) \subset L_s(\bG)$. Thus, the fractional operator is defined for all smooth and bounded functions.

\subsection{The Schwartz class of functions and the fractional operator}\label{s:Schwartz class}
Let
\[\mathcal{S}(\bG) = \{\phi\in C^\infty(\bG)\mid \phi\circ\exp \in \mathcal{S}(\algg)\}
\]
be the class of Schwartz  functions on $\bG$, see \cite[p. 1074]{KVW69}. Denote with $X_j$, resp. $Y_j$, $j=1,\dots, n$ the left-invariant, resp. right-invariant, vector fields,  see \eqref{e:inv vec fields}. Taking into account \cite[Proposition (1.26)]{FS82} any derivative with respect to  the exponential coordinates $x_j$, $j=1,\dots,n$, can be expressed as a linear combination with polynomial coefficients of the left- (or right-) invariant vector fields. Conversely, any  left- (or right-) invariant vector field can be written as a linear combination with polynomial coefficients of the  derivative $\frac {\partial}{\partial_{x_j}}$,  $j=1,\dots,n$.
In particular, we have the following alternative description, see \cite[Chapter I, D.]{FS82},
\begin{equation}\label{e:Schwartz class}
\mathcal{S}(\bG)=\{\phi\in \mathcal{C}^\infty(\bG) \mid P(x)X_I\phi \in L^\infty(\bG) \} \ = \ \{\phi\in \mathcal{C}^\infty(\bG) \mid P(x)Y_I\phi \in L^\infty(\bG) \}
\end{equation}
for any polynomial $P(x)$ and any multi-index $I$.  The space $\mathcal{S}(\bG)$ is a Fr\'echet space with topology defined via the family of norms, $k, N\in \mathbb{N}\cup \{0\}$
\[
\norm{\phi}_{N, k}=\sup_{|I|\leq k, g\in\bG} (1+|g|^2)^{N/2}|X_I\phi (g)|.
\]
For $u\in\mathcal{S}(\bG)$ the decay properties of $\frlap u(g)$ are given in the next proposition.

\begin{prop}\label{l:frac lapl def and decay}
If $u\in \mathcal{S}(\bG)$ then
 \begin{equation}\label{e:frac lapl decay for Schwartz}
{|g|^{Q+2s}}\vert\frlap u(g)\vert   \in L^{\infty}(\bG).
\end{equation}
In particular,  $\frlap u\in L^1(\bG)$.
If,  in  addition, $u$ has compact support, $\mathrm{supp} \, u \subset B_R$, then $\frlap u(g)$ has the following decay for large $|g|\geq 2R$,
\begin{equation}\label{e:frac lapl decay}
\vert\frlap u(g)\vert\leq \frac {C}{|g|^{Q+2s}}, \quad C= 2^{Q+2s+1}\sigma_Q R^Q\norm {u}_{L^{\infty}(\bG)}.
\end{equation}

\end{prop}

\begin{proof}

Let  $u\in \mathcal{S}(\bG)$. Let $4R=|g|\geq 4c$, where $c$ is the constant in \eqref{e:Taylor rem ord a}.  Splitting the integral \eqref{e:frac L via hypersing app app} over the ball $B(g,R)$ and its complement, then arguing as in \eqref{e:pv to non pv},  we have
\begin{equation}\label{e:Dsu split}
\frlap u(g)= 2\int_{\bG\setminus B(g,R)} \frac {u(g)-u(h)}{|g^{-1}h|^{Q+2s}}dh \ +\    \int_{B_R} \frac {2u(g)-u(gh)-u(gh^{-1})}{|h|^{Q+2s}}dh.
\end{equation}
The first integral is estimated by using the triangle inequality and \eqref{e:polar coord}, which give
\begin{multline}\label{e:decay on complemnt}
 \left| \int_{\bG\setminus B(g,R)} \frac {u(g)-u(h)}{|g^{-1}h|^{Q+2s}}dh\right|  \leq \left| u(g)\right|  \int_{\bG\setminus B(g,R)} \frac {1}{|g^{-1}h|^{Q+2s}}dh +  \int_{\bG\setminus B(g,R)} \frac {\left| u(h)\right|}{|g^{-1}h|^{Q+2s}}dh\\
 \leq  \frac{1 }{2s}\frac {\left| u(g)\right|}{R^{2s}} +  \frac {1}{R^{Q+2s}} \int_{\bG\setminus B(g,R)}  {\left| u(h)\right|}\,dh
= \frac{4^{Q+2s}\sigma_Q }{2s}\frac {\left| u(g)\right|}{|g|^{2s}} +  \frac{4^{Q+2s} \norm{u}_{L^1(\bG)}  } {|g|^{Q+2s}}\\
\leq \frac{C}{s}\left( \norm{u}_{Q,0}+\norm{u}_{Q+1,0}  \right) \frac{1 } {|g|^{Q+2s}}\leq \frac{C_Q}{s}\norm{u}_{Q+1,0}\frac{1 } {|g|^{Q+2s}}
\end{multline}
with a  constant $C_Q$.

To estimate the second integral in \eqref{e:Dsu split} we split the domain in two subsets and use the triangle inequality to obtain for some absolute constant $C$
\begin{multline*}
\left\vert\int_{B_R} \frac {2u(g)-u(gh)-u(gh^{-1})}{|h|^{Q+2s}}dh\right\vert \leq \left\vert\int_{B_1} \frac {2u(g)-u(gh)-u(gh^{-1})}{|h|^{Q+2s}}dh\right\vert \\
+ \left\vert\int_{B_R\setminus {B_1}} \frac {2u(g)-u(gh)-u(gh^{-1})}{|h|^{Q+2s}}dh\right\vert.
\end{multline*}
From  \eqref{e:2nd order difference} and  \eqref{e:2nd order difference remainder} for $|h|\leq 1$ we estimate the first integral, while we bound the second  by the triangle inequality and \eqref{e:polar coord}, to obtain
\begin{multline}\label{e:s to 0 on BR}
\left\vert\int_{B_R} \frac {2u(g)-u(gh)-u(gh^{-1})}{|h|^{Q+2s}}dh\right\vert
\leq  \sum_{i,j:\, d_i=d_j=1}\left\vert  (X_iX_ju)(g)\right\vert   \int_{ B_1} \frac {\left\vert x_i(h)x_j(h) \right\vert}{|h|^{Q+2s}}dh\\   +\  \int_{ B_1} \frac {\left\vert\mathcal{R}(h;g)\right\vert}{|h|^{Q+2s}}dh  + \left\vert\int_{B_R\setminus {B_1}} \frac {2u(g)-u(gh)-u(gh^{-1})}{|h|^{Q+2s}}dh\right\vert  \\
\leq \sum_{i,j:\, d_i=d_j=1}\left\vert  (X_iX_ju)(g)\right\vert \int_{ B_1}  \frac {|h|^{2}}{|h|^{Q+2s}}dh  \ + \  2T_2  \sum_{I:\, 1\leq |I|\leq 3, d(I)> 2} \sup_{|h'|\leq c|h|} |X_Iu(gh')|   \int_{ B_1} \frac {|h|^{d}}{|h|^{Q+2s}}dh\\
  + \int_{B_R\setminus {B_1}} \frac {\left\vert2u(g)-u(gh)-u(gh^{-1})\right\vert}{|h|^{Q+2s}}dh\\
\leq   C T_2\frac {\sigma_Q}{2 -2s}M_2(R)  \  + \ 4\sup_{h\in B_{R}(g)} \left | u (h) \right| \int_{\bG\setminus B_1} \frac {1}{|h|^{Q+2s}}dh
\leq   \frac{C_{Q}}{s(1-s)} \left(\sup_{h\in B_{R}(g)} \left | u (h) \right|   + M_{2}(R) \right),
\end{multline}where the constant $M_2(R)$ is defined, recalling $|g|=4R\geq 4c\geq 4$, by
\begin{equation}\label{e:M2 in Taylor}
M_{2}(R)=\sup_{B_{R}(g)} \{ \left |X_Iu \right| \  :  \   1\leq |I|\leq 3, \, 2\leq d(I)  \} \geq \sup_{B_{c}(g)} \{ \left |X_Iu \right| \  :  \   1\leq |I|\leq 3, \, 2< d(I)  \}
\end{equation}
and, further, we use the triangle inequality to see that for $|h'|\leq c|h|$ with $|h|\leq 1$ we have $|g^{-1}gh'|\leq c\leq R$ in order to estimate the supremum in the remainder \eqref{e:2nd order difference remainder} with $M_{2}(R)$, while for the last term we use $|g^{-1}gh|=|h|\leq R$.  Thus, we have
\begin{equation}\label{e:s to 0 on BR final}
\left\vert\int_{B_R} \frac {2u(g)-u(gh)-u(gh^{-1})}{|h|^{Q+2s}}dh\right\vert  \leq  \frac{C_{Q}}{s(1-s)} \frac {\norm{u}_{Q+2s,3}}{|g|^{Q+2s}}. \end{equation}

From \eqref{e:decay on complemnt} and \eqref{e:s to 0 on BR final} taking into account   \eqref{e:Dsu split}, the function $\frlap u$ has the claimed  decay at infinity for $u\in \mathcal{S}(\bG)$. The  fact that $\frlap u\in L^1(\bG)$ follows from \eqref{e:polar coord} and  the decay of $\frlap u$  shown above.

Suppose that $u$ is supported in the ball $B_{R_0}$, $ \mathrm{supp}\, u\subset B_{R_0}$. Hence,  for $|g|\geq 2R_0$ we have
\begin{multline}\label{e:frac lapl decay in proof}
\left\vert\int_{\bG\setminus B(g,\eps)} \frac {u(g)-u(h)}{|
h^{-1}\, g|^{Q+2s}}dh \right\vert = \left\vert \int_{B_{R_0}\cap {\bG\setminus B(g,\eps)}}\frac {u(h)}{|
h^{-1}\, g|^{Q+2s}}dh \right\vert\\
\leq \norm {u}_{L^{\infty}(\bG)}\int_{B_{R_0}}\frac {1}{|
h^{-1}\, g|^{Q+2s}}dh
\leq 2^{Q+2s}\norm {u}_{L^{\infty}(\bG)}|B_{R_0}|\frac {1}{|g|^{Q+2s}},
\end{multline}
after using the triangle inequality which gives for $|h|<{R_0}<2{R_0}<|g|$ the inequality
\[
| h^{-1} \, g|\geq |g|-|h|\geq |g|-{R_0}\geq |g|-\frac 12 |g|=\frac 12 |g|.
\]
Therefore, for $|g|\geq 2{R_0}$ we have
\[
\vert \frlap u(g)\vert \leq 2^{Q+2s+1}\norm {u}_{L^{\infty}(\bG)}|B_{R_0}|\frac {1}{|g|^{Q+2s}}.
\]

\end{proof}

We stress that for $u\in \mathcal{S}(\bG)$  it is not true that $\frlap u \in \mathcal{S}(\bG)$. In fact, even for a function in $C^\infty_0(\bG)$ one can only claim the decay in the above proposition. On the other hand, for  $\frlap u \in \mathcal{S}(\bG)$ we do have that  $\frlap u \in C^\infty(\bG)$ as shown in the next proposition.
 \begin{prop}\label{p:smoothness of frlap}
 If $ u \in \mathcal{S}(\bG)$ then $\frlap u \in \mathcal{C}^\infty (\bG)$.
 \end{prop}

 \begin{proof} As before, we use $X_j$, resp. $Y_j$, $j=1,\dots, n,$ to denote the left-invariant, resp. right-invariant, vector fields,  see \eqref{e:inv vec fields}.
 We recall \eqref {e:Schwartz class} regarding the equivalent definitions of the Schwartz class.  Thus, if $u\in \mathcal{S}(\bG)$ then  $Y_I u\in \mathcal{S}(\bG)$ for  any multi-index $I$. Furthermore, taking again into account \cite[Proposition (1.26)]{FS82}  we have that a function $u$ is $C^\infty(\bG)$-smooth (in the exponential coordinates)  if and only if $X_Iu$ (or $Y_Iu$)  are $C^\infty(\bG)$-smooth for all multi-indices $I$.

 Let us observe that the fractional operator commutes with the right-invariant vector fields $Y_j$. This is due to the fact that, if $R_h$ is the operator of right-translation defined by $(R_h u) (g)=u(gh)$, then $[Y_j,R_h]=0$. Therefore, if $u\in \mathcal{S}(\bG)$ then  for  any multi-index $I$  we have that $Y_I \frlap u$ exists and is given by
 \[
 Y_I \frlap u(g)=\frlap (Y_Iu)(g).
 \]
 This completes the proof.

 \end{proof}

\subsection{The carr\'e du champ operator}

Given a function $u\in \mathcal{D}^{s,2}(\bG)$, it will be convenient to define its fractional \emph{carr\'e du champ}
\begin{equation}\label{e:carre du champ}
{\Gamma_s (u)}(g)\overset{def}{=}  {\Gamma_s( u,u)}(g)\overset{def}{=}  \int_{\bG} \frac {|u(g)-u(h)|^2}{|h^{-1}\, g|^{Q+2s}}\, dh\geq 0
\end{equation}
and also
\begin{equation}\label{e:hor gradient}
{D_s u}(g)\overset{def}{=}\sqrt{\Gamma_s (u)(g)}.
\end{equation}
As usual, we will use the formula ${\Gamma_s (u,v)}$ for the symmetric bilinear form corresponding to ${\Gamma_s (u)}$,
\[
{\Gamma_s (u,v)}(g)\overset{def}{=} \int_{\bG} \frac {(u(g)-u(h))(v(g)-v(h))}{|h^{-1}\, g|^{Q+2s}}\, dh\
\]
Clearly, we have
\[
\norm {u}_{\mathcal{D}^{s,2}(\bG)}=\norm {D_su}_{L^2(\bG)}= \norm{{\Gamma_s (u)}}^{1/2}_{L^1(\bG)}  =\left( \int_\bG\int_{\bG} \frac {|u(g)-u(g)|^2}{|h^{-1}\, g|^{Q+2s}}\, dhdg\right )^{1/2}.
\]

\begin{lemma}\label{l:integral for Dirichlet form}
If $u,\, \phi\in C_0^\infty(\bG)$, then the function
\[
D(g,h)\overset{def}{=}\frac
{(u(g)-u(h))\left(\phi(g)-\phi(h)\right)}{| h^{-1}\, g|^{Q+2s}}\in L^1(\bG\times \bG).
\]
In particular,  ${\Gamma_s (u,v)}\in L^1(\bG)$.
\end{lemma}

\begin{proof}
 Suppose $R>0$ is such that $\operatorname{supp}\ u\subset B_{R}$.
The proof follows from Fubini's theorem once we show that
\[
\int_{\bG}\int_{\bG} |D(g,h)|\, dhdg = \int_{\bG\setminus B_{2R}}\int_{\bG} |D(g,h)|\, dhdg \ +\ \int_{B_{2R}}\int_{\bG} |D(g,h)|\, dhdg  <\infty.
\]
In order to shorten the formulas we will use the polarization formula, which allows a reduction to the case $u=\phi$. The arguments are very similar to the one used in \eqref{e:frac lapl decay in proof}.
For $|g|\geq 2R$, as in \eqref{e:frac lapl decay in proof} we have the bound
\begin{align}\label{e:D(g,h) estimate}
\int_{\bG\setminus B_{2R}}({D_s u})^2(g) \volg & =\int_{\bG\setminus B_{2R}}\int_{B_R} \frac {|u(h)|^2}{|h^{-1}\, g|^{Q+2s}}\, dhdg
\\
& \leq 2^{Q+2s}\norm {u}^2_{L^{\infty}(\bG)}|B_R|\int_{\bG\setminus B_{2R}}\frac {1}{|g|^{Q+2s}}\,dg <\infty.
\notag
\end{align}
Next, we consider the case $|g|<2R$. By the change of variable $h'^{-1}= h^{-1} g$, and the property $|h'^{-1}|=|h'|$, we have
\begin{equation}\label{e:small g case}
({D_s u})^2(g) = \int_{B_1} \frac {|u(gh)-u(g)|^2}{| h|^{Q+2s}} \, dh \ +\ \int_{\bG\setminus B_1} \frac {|u(gh)-u(g)|^2}{| h|^{Q+2s}}\, dh.
\end{equation}
We next show that each of the integrals in the right-hand side is uniformly bounded for $|g|<2R$. Indeed, for $|h|\leq 1$ the mean value inequality \eqref{e:MVT} gives
\[
|u(gh)-u(g)| \leq T_0|h|\sum_{j=1}^n\sup_{|h'|\leq 2R+c}|X_ju(h')|,
\]
taking into account that $d_i\geq 1$, $|h|\leq 1$, the triangle inequality, and the assumption $|g|<2R$. Since $|h|^{-Q+2(1-s)}\in L^1(B_1)$ the first integral in \eqref{e:small g case} is uniformly bounded for $|g|<2R$. On the other hand, due  to \eqref{e:polar coord} the integral over $|h|\geq 1$ can be estimated using that $|h|^{-(Q+2s)}\in L^1(\bG\setminus B_1)$,
 \begin{equation}\label{e:kernel integral exterior ball}
\int_{\bG\setminus B_1} \frac {1}{|h|^{Q+2s}} dh =\frac {\sigma_Q}{2s},
\end{equation}
while the numerator is bounded by $2\norm {u}^2_{L^\infty(\bG)}$.

\end{proof}

We can complete the

\begin{proof}[Proof of Theorem \ref{p:symm of d form and hypersing int}] The equivalence of the two definitions of the fractional operator is contained in Proposition \ref{l:frac lapl def and decay}.
Next, we show that the Dirichlet form is finite for functions with compact support and then prove \eqref{e:pv integral}.
Taking into account Lemma \ref{l:integral for Dirichlet form}, we can apply Fubini's theorem to the integral defining the Dirichlet form. After a change of variable, we  obtain
\[
\mathcal{D}(u,\phi){=}\int_{\bG}\int_{\bG} \frac
{(u(gh)-u(g))\left(\phi(gh)-\phi(g)\right)}{|h|^{Q+2s}}\,
{dh}dg= \int_{\bG}\frac {D(h)}{|h|^{Q+2s}}dh,
\]
where
\begin{multline*}
D(h)= \int_{\bG}
{(u(gh)-u(g))\left(\phi(gh)-\phi(g)\right)}\,
{dg}
 = -\int_{\bG}
{u(g)\left(\phi(gh)-\phi(g)\right)}\,
{dg} \\
+ \int_{\bG}
{u(gh)\left(\phi(gh)-\phi(g)\right)}\,
{dg}
=-\int_{\bG}
{u(g)\left(\phi(gh)-\phi(g)\right)}\,
{dg} + \int_{\bG}
{u(g)\left(\phi(g)-\phi(gh^{-1})\right)}\,
{dg}\\
=\int_{\bG} u(g)\left(2\phi(g)-\phi(gh)-\phi(gh^{-1}) \right)\, dg.
\end{multline*}
In the above, we have used that the product of $u$ and $\phi$, or their translates, is a function in $L^1(\bG)$. Therefore, since $\frlap \phi\in L^1(\bG)$ by Proposition \ref{l:frac lapl def and decay}, we obtain identity \eqref{e:pv integral}.

\end{proof}

We note explicitly some useful consequences of the above analysis.
First,  from \eqref{e:D(g,h) estimate} we have that if $\phi\in C^{\infty}_0 (\bG)$ with $\mathrm{supp}\, \phi\subset B_R$, then the carr\'e du champ function  \eqref{e:hor gradient} has the following decay at infinity
\begin{equation}\label{e:grad square decay}
D_s \phi (g)=\sqrt{\Gamma_s (\phi)(g)}  \leq 2^{(Q+2s)/2}|B_R|^{1/2}\frac {\norm {\phi}_{L^{\infty}(\bG)}}{|g|^{(Q+2s)/2}}.
\end{equation}
Furthermore,  for $u, v\in C^{\infty}_0 (\bG)$, we have
\begin{equation}\label{e:product rule}
\frlap (uv)=u\frlap v + v\frlap u - 2\Gamma_s(u,v)
\end{equation}
and
\begin{equation}\label{e:int by parts}
\int_{\bG}  \Gamma_s(u,v) \, dg= \mathcal{D}(u,v)=\int_{\bG}v(g)\, \frlap u (g) dg=\int_{\bG}u(g)\, \frlap v (g) dg.
\end{equation}
Finally, we observe that the Sobolev embedding theorem \ref{t:frac homog Sobolev and Rellich embed} can be trivially stated  in the equivalent  form
\begin{equation}\label{e:fracsobolev as SGN}
\left ( \int_{\bG} |u|^{2^*(s)}(g){dg}\right)^{1/2^*(s)}\leq S\ \norm{D_s u}_{L^2(\bG)}.
\end{equation}

\subsection{Group invariances of the fractional operator}\label{ss:invariances frac operator}

As observed before, the fractional operator commutes with left-translations, but not necessarily with right-translations. In fact, if we let $L_{g_0} g= g_0\, g$, then in view of the  following identity,
\[
\frlap (L_{g_0}u)(g)=\int_{\bG} \frac {2u(g_0g)-u(g_0gh)-u(g_0gh^{-1})}{|h|^{Q+2s}}dh=L_{g_0}(\frlap u)(g),
\]
we have
\begin{equation}\label{frac translation prop}
L_{g_0}(\frlap u)=\frlap (L_{g_0}u).
\end{equation}
Also, it is easy to show that
\begin{equation}\label{frac dilation prop}
\frlap(\delta_\lambda u)=\lambda^{2s}\delta_\lambda(\frlap u).
\end{equation}

\subsection{ The limits of the fractional operator at 0 and 1. }\label{ss:limits} In this section we prove Theorem \ref{p:limits frac operator}.
As a preparation, we start with some examples. Recall that the  homogeneous norm defined in \cite[Theorem 2]{HS90} is invariant under orthogonal transformations in each of the vector spaces spanned by the  homogeneous vectors $X_i$ in $\algg$ of the same homogeneity, see also  \cite[p.1074]{KVW69} and  \cite[p. 604]{KV71}. In fact,  there exists an $\varepsilon>0$, such that, for any $r<\varepsilon$ the following formula defines a homogeneous norm
\begin{equation}\label{e:HS homog norm}
|g|=|g|_r\overset{def}{=}\inf_{\lambda>0}\{\lambda \mid \, \norm{\delta_{1/\lambda} x} <r \},
\end{equation}
where $\exp g =x$ and  $\norm{x} =\left(\sum_{i=1}^n x_i^2 \right)^{1/2}$. This norm satisfies the triangle inequality, cf. \eqref{e:homog norm} and \eqref{e:triangle ineqs}. Furthermore, the unit gauge ball with respect to the norm $|g|_r$ is the Euclidean ball of radius $r$,  $\{x\mid \, \norm {x}<r\}$, in the fixed homogeneous coordinates. Thus, we could assume, when necessary, that the fixed  homogeneous norm $|g|$ is a poly-radial function of the coordinates \eqref{e:homog coord} split according to their homogeneity.

In the Abelian case, i.e., $R^n$ with the addition and dilations given by $$\delta_\lambda(x_1,\dots,x_n)=(\lambda^{d_1}x_1,\dots, \lambda^{d_1}x_n),$$ such a norm was written explicitly in \cite[Remark 1]{FR66}. As another example, consider the Heisenberg group $\mathbb H^n = (\mathbb{C}^n \times \R, \circ)$, with group law
\begin{equation}\label{e:H-type Iwasawa groups}
  (z_0, t_0)\circ(z, t)\ =\ (z_0 + z, t_0+ t + \frac12\ \mathrm{Im}\,  z_0\cdot \bar z).
\end{equation}
where $z$, $z_0\in \mathbb{C}^n$ and $t, t_0\in \mathbb{R}$ with $z\cdot z' = \sum_{j=1}^n z_j {z}'_j$. The exponential map is the identity and  the parabolic dilations are defined by
$\delta_\lambda(z,t)=(\lambda z, \lambda^2 t)$. The corresponding homogeneous dimension is $Q= 2n + 2$. Identifying  $\mathbb{C}^n \times \R\cong \R^{2n} \times \R$,  and $z=x+iy \in \mathbb{C}^n$ with $(x,y)\in \mathbb{R}^{2n}$, a real basis for the Lie algebra of right-invariant vector fields on $\mathbb H^n$ is given by
\begin{equation}\label{e:basis for Hn}
Y_j =\frac{\partial}{\partial x_j} +\frac {y_j}{2} \frac{\partial}{\partial t}, \quad Y_{n+j}=  \frac{\partial}{\partial y_j} -\frac { x_j}{2} \frac{\partial}{\partial t},\quad Y_{2n+1}=\frac{\partial}{\partial t}, \quad j=1,...,n,
\end{equation}
while the left-invariant vector fields are
\[
X_j=\frac{\partial}{\partial x_j} -\frac {y_j}{2} \frac{\partial}{\partial t}, \quad X_{n+j}\equiv Y_j =  \frac{\partial}{\partial y_j} +\frac {x_j}{2} \frac{\partial}{\partial t},\quad X_{2n+1}=\frac{\partial}{\partial t}, \quad j=1,...,n.
\]
The first $2n$ vectors are of homogeneity one and the last vector is homogeneous of degree 2.
In this case, following \eqref{e:HS homog norm}, we have
\[
|(z,t)|_r = \frac {1}{\sqrt 2}\left[  \frac {|z|^2}{r^2} +\sqrt{\frac {|z|^4}{r^4}+4\frac {t^2}{r^2}}\right]^{1/2}=\lambda_r,
\]
where $\lambda_r$ is the positive root of the equation $\lambda^4r^2-\lambda^2|z|^2-t^2=0$.
Notice that $|(z,t)|_r<1$ is equivalent to $ |z|^2+t^2<r^2$.
One can also use the Kor\'{a}nyi gauge
\[
|(z,t)|=\left[|z|^4+ 16t^2 \right]^{1/4}.
\]
As mentioned in the introduction, with this choice of the gauge, we obtain the geometric (pseudo-hermitian conformally invariant) fractional sub-Laplacian, see  \cite{FGMT} and  \cite[Proposition 4.1]{RTaim}.

Gauges with similar symmetries have been considered in the setting of Carnot groups in \cite{Barb11} and \cite{CMSV21}. In particular, these papers  contain results in the spirit of \cite{BBM01} in the case of fractional Sobolev spaces  \cite[Theorem 3.6]{Barb11}  and Orlicz spaces \cite[Theorem 1.1]{CMSV21}.

We  turn to the proof of Theorem \ref{p:limits frac operator}. Let us mention that in order to compute the limits of the fractional operator as $s\rightarrow0+$ and $s\rightarrow 1-$ we only use that the gauge is invariant under a smaller group of  symmetries than the assumption in Theorem \ref{p:limits frac operator},  see Remark \ref{r:gauge smaller group symm}.
\begin{proof} [ Proof of Theorem \ref{p:limits frac operator}]
 Since we assume that the homogeneous norm $|g|$ is a radial function of the coordinates of homogeneity one,  if $x_i(h)$ and $x_j(h)$ are distinct among the coordinates that are homogeneous of order one, we have for $i\not=j$
\begin{equation}\label{e:xixj integrals}
\int_{B_1} \frac {x_i(h)x_j(h)}{|h|^{Q+2s}}dh=0\quad\text{and}\quad \int_{B_1} \frac {x^2_i(h)}{|h|^{Q+2s}}dh=\frac {\tau_m}{2m(1-s)},
\end{equation}
where $\tau_m$ is the constant defined in \eqref{e:tau_m def}.  Indeed, for a fixed $i$, a change of variables $y_j=x_j$ when $j\not= i$ and $y_i=-x_i$, taking into account the assumed invariance of the norm and the fact that the Haar measure is the Lebesgue measure gives
\[
\int_{B_1} \frac {x_i(h)x_j(h)}{|h|^{Q+2s}}dh=-\int_{B_1} \frac {x_i(h)x_j(h)}{|h|^{Q+2s}}dh,
\]
while the polar coordinates formula \eqref{e:polar coordinates integral}  and the homogeneity of $x_i$ allows us to compute the integral in the second formula,
\begin{multline*}
\int_{B_1} \frac {x^2_i(h)}{|h|^{Q+2s}}dh=\int_0^1\int_{|h|=1}\frac {r^{Q-1}r^2 x_i^2(h)}{r^{Q+2s}}\,d\sigma(h)dr=\int_0^1\frac {dr}{r^{2s-1}}\, \int_{|h|=1} x_i^2(h)\,d\sigma(h)\\
=\frac {1}{2(1-s)}\frac {1}{m}  \int_{|h|=1} \sum_{j=1}^m x_j^2(h)\,d\sigma(h).
\end{multline*}

Let $u\in \mathcal{ S}(\bG)$. First, we will compute the limit of $\frlap u(g)$ when $s\rightarrow 1-$.  The estimates are similar to the ones in proposition \ref{p:defs of frac operator}, except now, due to the extra assumption on the norm, we can compute precisely the integrals of the terms in \eqref{e:2nd difference quot}  involving $|h|^2$. With the help of \eqref{e:2nd order difference} and  \eqref{e:xixj integrals} we compute
\begin{multline}\label{e:Dsu formula}
\frlap u(g)= \int_{B_1} \frac {2u(g)-u(gh)-u(gh^{-1})}{|h|^{Q+2s}}dh+ \int_{\bG\setminus B_1} \frac {2u(g)-u(gh)-u(gh^{-1})}{|h|^{Q+2s}}dh \\
=-\frac {\tau_m}{2m(1-s)}\sum_{i=1}^m X_i^2u(g) \ +\  \int_{ B_1} \frac {\mathcal{R}(h;g)}{|h|^{Q+2s}}dh
+ \int_{\bG\setminus B_1} \frac {2u(g)-u(gh)-u(gh^{-1})}{|h|^{Q+2s}}dh.
\end{multline}
The last two terms in the above formula for the fractional operator $\frlap u(g)$ are bounded as $s\rightarrow 1-$. Indeed, for $d >2$ and \eqref{e:polar coord} we obtain
\[
\int_{B_1} \frac {|h|^{d}}{|h|^{Q+2s}}dh=\frac {\sigma_Q}{d-2s}.
\]
Therefore, using also \eqref{e:2nd order difference remainder}, \eqref{e:kernel integral exterior ball} and the triangle inequality, we have
\begin{equation}\label{e:Dsu formula remainder}
\int_{ B_1} \frac {|\mathcal{R}(h;g)|}{|h|^{Q+2s}}dh
+ \int_{\bG\setminus B_1} \frac {|2u(g)-u(gh)-u(gh^{-1})|}{|h|^{Q+2s}}dh\leq 2T_2M_{2}(c)\frac {\sigma_Q}{d -2s}+ 4\norm{u}_{L^\infty }\frac {\sigma_Q}{2s},
\end{equation}
where  $c\geq 1$ by \cite[p. 29]{FS82}  is the constant appearing in the formula for the Taylor remainder \eqref{e:Taylor rem ord 2} and the constant  $M_{2}(c)$ is defined in \eqref{e:M2 in Taylor}.
This completes the proof of the second identity in \eqref{e:limits frac operator}.

We turn to the case  $s\rightarrow 0$. Suppose $ \mathrm{supp}\, u \subset B_{R_0}$. For $g\in \bG$  let $R=R_0+|g|+c>1$.  As before, we estimate the integral defining the fractional operator by splitting it into an integral over the ball $B_R$ and its complement,
\begin{equation}\label{e:Dsu formula s at 0}
\frlap u(g)= \int_{B_R} \frac {2u(g)-u(gh)-u(gh^{-1})}{|h|^{Q+2s}}dh\ +\  \int_{\bG\setminus B_R} \frac {2u(g)-u(gh)-u(gh^{-1})}{|h|^{Q+2s}}dh
\end{equation}
Taking into account that  $ \mathrm{supp}\, u \subset B_{R_0}$, the triangle inequality, and the definition of $R$ we have, recall \eqref{e:polar coord} and \eqref{e:Dsu formula remainder},
\begin{equation}\label{e:s to 0 on complement of BR}
 \int_{\bG\setminus B_R} \frac {2u(g)-u(gh)-u(gh^{-1})}{|h|^{Q+2s}}dh = 2u(g)\int_{\bG\setminus B_R} \frac {1}{|h|^{Q+2s}}dh= 2u(g)\frac {\sigma_Q R^{-2s}}{2s}.
\end{equation}
When multiplied by $s/\sigma_Q$, the last quantity converges to $u(g)$ as $s\rightarrow 0+$, which equals the limit claimed in Theorem \ref{p:limits frac operator}.

Next we will show that the first integral in \eqref{e:Dsu formula s at 0}  multiplied by $s$ converges to 0 as $s\rightarrow 0+$. Splitting the domain in two subsets, from the triangle inequality we have the estimate
\begin{multline*}
\left\vert\int_{B_R} \frac {2u(g)-u(gh)-u(gh^{-1})}{|h|^{Q+2s}}dh\right\vert \leq \left\vert\int_{B_1} \frac {2u(g)-u(gh)-u(gh^{-1})}{|h|^{Q+2s}}dh\right\vert \\
+ \left\vert\int_{B_R\setminus {B_1}} \frac {2u(g)-u(gh)-u(gh^{-1})}{|h|^{Q+2s}}dh\right\vert.
\end{multline*}
Using \eqref{e:2nd order difference remainder} for $|h|\leq 1$ and the trivial inequality $1\leq |h|^d$ otherwise,  taking into account also \eqref{e:Dsu formula}, the following bounds follow
\begin{multline}\label{e:s to 0 on BR}
\left\vert\int_{B_R} \frac {2u(g)-u(gh)-u(gh^{-1})}{|h|^{Q+2s}}dh\right\vert \leq \\
\leq \left\vert\frac {\tau_m}{2m(1-s)}\sum_{i=1}^m X_i^2u(g) \right\vert \ +\  \int_{ B_1} \frac {\left\vert\mathcal{R}(h;g)\right\vert}{|h|^{Q+2s}}dh  + \left\vert\int_{B_R\setminus {B_1}} \frac {2u(g)-u(gh)-u(gh^{-1})}{|h|^{Q+2s}}dh\right\vert  \\
\leq  \left\vert\frac {\tau_m}{2m(1-s)}\sum_{i=1}^m X_i^2u(g) \right\vert \ + \ 2 T_2M_{2}\frac {\sigma_Q}{d -2s} \  + \ 4\norm{u}_{L^\infty } \int_{B_R} \frac {|h|^d}{|h|^{Q+2s}}dh\\
= \left\vert\frac {\tau_m}{2m(1-s)}\sum_{i=1}^m X_i^2u(g) \right\vert \  + \left(4\norm{u}_{L^\infty }+2 T_2M_{2} \right) \sigma_Q\frac {R^{d-2s}}{d-2s}.
\end{multline}
When multiplied by $s$, the limit  of the last line when  $s\rightarrow 0+$ is 0.
Therefore,   the first identity in \eqref{e:limits frac operator} is proved, as well.

\end{proof}

\begin{rmrk}\label{r:gauge smaller group symm}
The proof of Proposition \ref{p:limits frac operator} shows that we can weaken the radial symmetry with respect to the coordinates of homogeneity one by requiring that the norm is an even function of these coordinates. However, in this case the second formula in \eqref{e:limits frac operator} becomes
\[
\lim_{s\rightarrow 1-}  {2m(1-s)}\frlap u(g)=-\sum_{i=1}^m \hat{\tau}^2_{i}X_i^2u(g), \qquad \hat{\tau}^2_{i}=\int_{|h|=1} x_i^2\, d\sigma.
\]
In this case, it could be useful to re-scale the originally fixed vector fields (of homogeneity one) $X_1, \dots, X_m$  or consider the operator
\[
\hat{\frlap } u(x)=\frlap v(x), \qquad v(x)=u(\frac {x_1}{{\hat\tau}_{1}}, \dots, \frac {x_m}{{\hat\tau}_{m}}, x_{m+1},\dots,x_n).
\]
\end{rmrk}

\begin{rmrk}
We could consider the fractional operator $\alpha_{Q,m,s}\frlap$ instead with
\begin{equation}\label{e:definition of alpha}
\alpha_{Q,m,s}=-\frac { 2m  (1-s) s^2}{\tau_m} -\frac {s (1-s)^2}{\sigma_Q}=-s(1-s)\left( \frac {2 m  s}{\tau_m} +\frac { (1-s)}{\sigma_Q}\right).
\end{equation}
Although defined in somewhat artificial manner, this definition is motivated by several examples, so that in the limits when $s$ goes to 0 or $1-$ we obtain, respectively, the identity or the sum of squares of the vector fields of degree 1.
\end{rmrk}

With the help of  Theorem \ref{p:limits frac operator} we prove the following connection between the fractional Dirichlet norm and the $L^2$ Dirichlet norm of a smooth function with compact support.

\begin{prop}\label{t:BBM and MS} Let  $m$ be the number of left-invariant vector fields $X_j$ of homogeneity one. Suppose the fixed  homogeneous norm $|g|$ is a radial function of the coordinates \eqref{e:homog coord} of homogeneity one.
For $u\in C^\infty_0(\bG)$ we have
\begin{equation}\label{e:BBM for Cinfty0}
\lim_{s\rightarrow 1-} \frac {2m(1-s)}{\tau_m}\norm{u}^2_{\mathcal{D}^{s,2}(\bG)}=\int_{\bG}\sum_{i=1}^m |X_iu(g)|^2 dg.
\end{equation}
and
\begin{equation}\label{e:MS for Cinfty0}
\lim_{s\rightarrow 0+} \frac {s}{\sigma_Q}\norm{u}^2_{\mathcal{D}^{s,2}(\bG)}=\int_{\bG} |u(g)|^2 dg.
\end{equation}
\end{prop}

\begin{proof} We recall \eqref{e:def frac norm}. In order to prove \eqref{e:BBM for Cinfty0}, we use \eqref{e:pv integral} and \eqref{e:limits frac operator} which imply
\begin{multline}
\lim_{s\rightarrow 1-} \frac {2m(1-s)}{\tau_m}\mathcal{D}_s(u,u)=\lim_{s\rightarrow 1-} \frac {2m(1-s)}{\tau_m}\int_{\bG}u(g)\, \frlap u (g) dg\\
=- \int_{\bG}u(g)\, \sum_{i=1}^m X_i^2u(g) dg =\int_{\bG}\sum_{i=1}^m |X_iu(g)|^2 dg
\end{multline}
since $X_i^*=-X_i$ by \eqref{e:skew-adjoint v.f.}. To justify the limit, by Lebesgue's dominated convergence theorem it is enough to show that $ {(1-s)}|u(g)\, \frlap u (g)|$ is bounded by an $L^1(\bG)$ function, which is independent of $s$. For $u\in C^\infty_0(\bG)$, this fact follows from  \eqref{e:Dsu formula} and \eqref{e:Dsu formula remainder}.

The identity \eqref{e:MS for Cinfty0} follows similarly from Lebesgue's dominated convergence theorem using  \eqref{e:s to 0 on BR} and \eqref{e:s to 0 on complement of BR}.
\end{proof}

It should be noted that in the case of Carnot groups, the sum of squares of vector fields of homogeneity one that appears in Theorem \ref{p:limits frac operator}  is the left-invariant sub-Laplacian, while the integral in \eqref{e:BBM for Cinfty0} is the (square of the) $L^2$ norm of the horizontal gradient. In particular,  part of proposition \ref{t:BBM and MS}  generalizes the known Maz’ya and Shaposhnikova's results on Carnot group mentioned below. In the general setting considered here, the sum of squares of vector fields of homogeneity one is not a hypoelliptic operator and a  heat semigroup is not well-defined.

In the case of Carnot groups, results similar to proposition \ref{t:BBM and MS}  were obtained in \cite[Section 4 and Theorem 1.1] {CMSV21} and \cite[Theorem 3.6]{Barb11}.  More recently, again in the framework of Carnot groups, universal characterizations of BV and Sobolev spaces (with dimension independent constants) were obtained using heat semigroup techniques in \cite{GT24, GT23CAG}, and a dimensionless generalisation of the Maz’ya and Shaposhnikova's results was proven in \cite[Theorem 1.5]{GT24}.    In a different setting, in the works \cite{GT0NA} and \cite{BGT22} the authors considered the non-symmetric semigroup $P^\sA_t$ associated with the (possibly) degenerate evolution
$$\mathscr A u  - \de_t u \overset{def}{=} \operatorname{tr}(Q \nabla^2 u) + \langle BX,\nabla u\rangle - \de_t u,$$
where $Q, B\in \mathbb M_{n\times n}(\R)$, $Q\ge 0$ and $Q=Q^\star$, and $Q, B$ satisfy the Kalman assumption for complete controllability: $\operatorname{Ker} Q$ does not contain any nontrivial invariant subspace of $B^\star$, see \cite{Ho}. The authors introduced
the Besov space $\mathfrak B^\sA_{s,p}$  of all functions $f\in L^p(\Rn)$ such that
\begin{equation}\label{sn}
\mathscr N^\sA_{s,p}(f) = \left(\int_0^\infty  \frac{1}{t^{\frac{s p}2 +1}} \int_{\Rn} P^\sA_t\left(|f - f(x)|^p\right)(x) dx dt\right)^{\frac 1p} < \infty.
\end{equation}
Note that, in the non-degenerate setting of the classical heat semigroup, $Q = I_n, B = O_n$, one obtains
\[
\mathscr N^\Delta_{s,p}(f) = \frac{2^{sp}\Gamma(\frac{n+p}2)}{\pi^{\frac n2}}\int_{\Rn}\int_{\Rn} \frac{|f(x)-f(y)|^p}{|x-y|^{n+sp}} dx dy,
\]
which shows that for $p=2$ the seminorm \eqref{sn} recovers \eqref{e:def frac norm}. In this respect, see the semigroup based theory of nonlocal isoperimetric inequalities developed in \cite{GTiso}. The main result in \cite[Theorem 1.1]{BGT22} states that, when $\operatorname{tr} B \ge 0$, if $f\in \underset{0<s<1}{\bigcup}\mathfrak B^\sA_{s,p}$ the following dimensionless version of Maz'ya \& Shaposhnikova result holds:
\begin{equation}\label{ms}
\underset{s\to 0^+}{\lim} s \mathscr N^\sA_{s,p}(f)^p = \begin{cases}
\frac{4}{p} ||f||_p^p,\quad \mbox{if }  \operatorname {tr}B = 0,
\\
\frac{2}{p}||f||_p^p,\quad\ \mbox{if } \operatorname {tr} B >0.
\end{cases}
\end{equation}

\subsection{The Dirichlet space and proof of Theorem \ref{t:frac homog Sobolev and Rellich embed} a)}\label{ss:frac homog Sobolev and Rellich embed}
The Dirichlet space $\mathcal{D}^{s,2}(\bG) $ associated to the quadratic form $\mathcal{D}$ given in \eqref{e: dirichlet form} is defined as the completion of $C_0^\infty(\bG)$ with respect to the norm \eqref{e:def frac norm}. The space $\mathcal{D}^{s,2}(\bG) $ is a space of functions in $L^{2^*(s)}(\bG)$ due to the fractional Sobolev inequality Theorem \ref{e:fracsobolev} a).
For completeness, we include the short  proof of the inequality  based on the well-known interpolation argument of Hedberg and Stein, see \cite{He72}. We follow \cite[Proposition 15.5]{Po16}, where the inequality is proven in the Euclidean setting, but the argument is easily adaptable to homogeneous groups. We note that we do not assume that the homogeneous norm $|g|$ is a radial function of the coordinates of homogeneity one.

\begin{proof}[Proof of Theorem \ref{t:frac homog Sobolev and Rellich embed} a)]
For $g\in\bG$ and any $r>0$  after several applications of H\"older's inequality and the trivial inequality $$|g^{-1}h|^{Q+2s}\leq r^{Q+2s} \qquad \text{for}\quad |g^{-1}h|\leq r,$$ we obtain  the following estimate of $|u(g)|$
\begin{multline}
|u(g)|\leq \fint_{B(g,r)} |u(g)-u(h)|\,dh   +  \fint_{B(g,r)} |u(h)|\,dh\\
\leq \left( \fint_{B(g,r)} |u(g)-u(h)|^2\,dh \right)^{1/2}   +  \left( \fint_{B(g,r)} |u(h)|^{2^*(s)}\,dh\right)^{1/2^*(s)}
\leq Ar^a +Br^{-b},
\end{multline}
letting
\[
 a=s, \quad b=\frac {Q}{2^*(s)}=\frac {Q-2s}{2},
\]
 and
\[
 A=\left( \frac {1}{\sigma_Q}\int_{\bG} \frac{|u(g)-u(h)|^2}{|
h^{-1}\, g|^{Q+2s}}\,dh \right)^{1/2}\quad \text{and}\quad B=\left( \frac {1}{\sigma_Q} \int_{\bG} |u(h)|^{2^*(s)}\,dh\right)^{1/2^*(s)}.
 \]
From Lemma \ref{l:simple lemma} and the trivial identity
\[
\frac {1}{2^*(s)}+\frac sQ=\frac 12
\]
we obtain the inequality
\[
|u(g)|^{2^*(s)}\leq C^{{2^*(s)}}_{Q,s} \left( \int_{\bG} \frac{|u(g)-u(h)|^2}{|
h^{-1}\, g|^{Q+2s}}\,dh \right)\, \left(  \int_{\bG} |u(h)|^{2^*(s)}\,dh\right)^{1-\frac {2}{2^*(s)}},
\]
where
\[
C_{Q,s} = {\sigma_Q^{-\frac 14\frac{Q+2s}{Q-2s} }}\left[ \left(\frac ba\right)^ {\frac {a}{a+b}}  +  \left(\frac ba\right)^ {-\frac {b}{a+b}} \right]={\sigma_Q^{-\frac 14\frac{Q+2s}{Q-2s} }}  \left[ \left(\frac {Q-2s}{2s}\right)^ {\frac {2s}{Q}}  +  \left(\frac {Q-2s}{2s}\right)^ {-\frac {Q-2s}{4}} \right].
\]
Therefore,
\[
\left ( \int_{\bG} |u|^{2^*(s)}(g){dg}\right)^{1/2^*(s)}\leq C^{{Q}/(Q-2s)}_{Q,s}
\left(\int_{\bG}\int_{\bG} \frac {|u(g)-u(h)|^2}{|
h^{-1}\, g|^{Q+2s}}\, {dg}dh \right )^{1/2}.
\]
\end{proof}

\begin{lemma}\label{l:simple lemma}
Let $a$, $b$, $A$ and $B$ be positive constants. If $f(r)=Ar^a+Br^{-b}$,  $r>0$, then
\[
f(r)\geq \left[ \left(\frac ba\right)^ {\frac {a}{a+b}}  +  \left(\frac ba\right)^ {-\frac {b}{a+b}} \right]\, A^{\frac {b}{a+b}} B^ {\frac {a}{a+b}}.
\]
Furthermore, the right hand-side is the minimum of $f$ on $r>0$.
\end{lemma}

\begin{proof}
It is enough to prove the inequality for $A=B=1$ since simple scaling gives the general case. The only critical point of $f$ is its minimum and it is achieved at
$$r_0=\left(\frac ba\right)^{\frac {1}{a+b}}.$$
The value $f(r_0)$ gives the claimed bound.
\end{proof}

\begin{rmrk}\label{r:uniform frac sobolev}
One might wonder if we can get the full Bourgain et al. \cite{BBM01} and Maz'ya \& Shaposhnikova \cite{MS02} type results.
 Thus, it is reasonable to ask if a version of the following uniform with respect to $s\in(0,1)$ fractional Sobolev inequality in Euclidean space \cite[Theorem 1]{MS02}
\[
\norm{u}^2_{L^{2^*(s)}(\R^n)}\leq C(n,2) \frac {s(1-s)}{n-2s}\norm{u}^2_{{\mathcal{D}^{s,2}(\Rn)}}
\]
holds true for the fractional Sobolev inequality on homogeneous groups considered in this paper?  The answer in general is negative, as the following example shows.
\end{rmrk}

\begin{exm} \label{ex:non carnot}
Consider $\bG=\R^2$ with its Abelian structure, but with the parabolic dilations
$\delta_\lambda(x,y)=(\lambda x, \lambda^2 y)\in \R\times \R$. A homogenous norm is given by $|(x,y)|=|x|+|y|^{1/2}$. The vector field $X=\partial_x$ is of homogeneity one, while $Y=\partial_y$ is of homogeneity two. In particular, the fractional operator  is given by
\[
\frlap u (x,y)= p.v. \int_{\R^2} \frac {u(x,y)-u(\xi,\eta)}{\left (|x-\xi |+|y-\eta |^{1/2} \right)^{3+2s}}\, d\xi d\eta.
\]
By theorem \ref{t:BBM and MS} we have  $ \lim_{s\rightarrow 1-} \frac {2(1-s)}{\tau_1}\frlap u(x,y)=-u_{xx}(x,y)$. The fractional Sobolev inequality takes the form
\[
\left ( \int_{\R^2} |u (g)|^{6/(3-2s)}{dg}\right)^{(3-2s)/6}\leq S
\left(\int_{\R^2}\int_{\R^2} \frac {|u(g)-u(h)|^2}{\left (|x-\xi|+|y-\eta|^{1/2} \right)^{3+2s}}\, {dg}dh \right )^{1/2}
\]
with $g=(x,y)$ and $h=(\xi,\eta)$. Notice that $S$ depends on $s$.
If we were to have inequality
\[
\left ( \int_{\R^2} |u(x,y)|^{6/(3-2s)} {\, dg}\right)^{(3-2s)/6}\leq C(1-s)^{1/2}
\left(\int_{\R^2}\int_{\R^2} \frac {|u(x,y)-u(\xi,\eta)|^2}{\left (|x-\xi|+|y-\eta|^{1/2} \right)^{3+2s}}\, {dg}dh \right )^{1/2}
\]
with a constant $C$ independent of $s$ when $s\rightarrow 1-$, as in the Euclidean case of remark \ref{r:uniform frac sobolev}, then  from \eqref{e:BBM for Cinfty0} taking the limit $s\rightarrow 1-$ would result in
\[
\left ( \int_{\R^2} |u(x,y)|^6{dxdy}\right)^{1/6}\\
\leq C_1 \left(\int_{\R^2}|\partial_x u (x,y)|^2 \, dxdy \right )^{1/2}.
\]
The latter inequality does not hold if we take $u=\phi(x)\psi(y)$, since then we should have
\begin{multline}
\left ( \int_{\R} |\phi(x)|^6{dx}\right)^{1/6}  \left ( \int_{\R} |\psi(y)|^6{dy}\right)^{1/6}=
\left ( \int_{\R^2} |u(x,y)|^6{dxdy}\right)^{1/6}\\
\leq C_1
\left(\int_{\R^2}|\partial_x u (x,y)|^2 \, dxdy \right )^{1/2}=C_1\left ( \int_{\R} |\partial_x\phi(x)|^2{dx}\right)^{1/2}\left ( \int_{\R} |\psi(y)|^2{dy}\right)^{1/2},
\end{multline}
which is false.
\end{exm}

In the following lemma we show the existence of a useful truncation function. In particular, the associated multiplication operator has  properties which will be used in Theorem \ref{t:trancation thrm}.
\begin{lemma}\label{l:cutof for trancation}
For $R>0$ there exists $\phi\equiv \phi_R\in C^\infty(\bG)$, such that,
\begin{equation}\label{e:cutof for trancation}
\phi\equiv 0 \quad\text{on}\quad B_R, \quad \phi\equiv 1 \quad\text{on}\quad \bG\setminus {B_{2R}}, \quad \text{and} \quad
|\phi(g)-\phi(h)|\leq \frac {2}{R}|h^{-1}g|.
\end{equation}
\end{lemma}

\begin{proof}
The construction is very similar to that  of the bump function in \cite[Section 3.2]{GLV22}. Let $\eta$ be a smooth non-decreasing function on the real line with $\mathrm{supp}\, \phi\subset [R,+\infty)$, $\eta(t)=1$ for $t\geq 2R$,  and $|\eta'(t)|\leq 2/R$, so that,
\[
|\eta(t_1)-\eta(t_2)|\leq \frac {2}{R}|t_1-t_2|.
\]
Hence, the function $\phi(x)=\eta(|x|)$ satisfies all the claimed conditions, noting that
\[
|\phi(g)-\phi(h)|=  |\eta(|g|)-\eta(|h|)|  \leq \frac {2}{R}\left\vert \,|g|-|h|\, \right\vert \leq  \frac {2}{R}|h^{-1}g|,
\]
where the last inequality follows by the triangle inequality \eqref{e:triangle ineqs}.

\end{proof}

\begin{thrm}\label{t:trancation thrm}
The space $\mathcal{D}^{s,2}(\bG)$ can be identified with the space of functions $u\in L^{2^*(s)}(\bG)$ with finite $[\, .\, ]_{s,2}$ seminorm.
\end{thrm}

\begin{proof}
The proof can be done exactly as in the Euclidean case, using  the existence of the truncation function  in Lemma \ref{l:cutof for trancation}, see \cite[Theorem 2.2]{dPQ20} or \cite[Theorem 3.1]{BG-CV21}.
We sketch a slight variation relying on Proposition \ref{l:conv and Sobolev norm}. A calculation shows that for $\phi_R$  as in Lemma \ref{l:cutof for trancation}, there exists a constant $C=C(Q,s)$ such that for $u\in L^{2^*(s)}(\bG)$ with $[u]_{s,2}<\infty$, we have
\begin{equation}\label{e:trancation lemma}
[\phi_R u]_{s,2}\leq C \left[  \int_{\bG\setminus B_{2R}}\int_{\bG\setminus B_{2R}} \frac {|u(g)-u(h)|^2}{|h^{-1}g|^{Q+2s}}\volh\volg \ +\
\left(\int_{\bG\setminus B_{2R}}|u(g)|^{2^*(s)}\volg\right) ^{\frac {2}{2^*(s)}}
\right].
\end{equation}
Therefore, we have $\lim_{R\rightarrow \infty}[\phi_R u]_{s,2}=0$. On the other hand, a calculation shows that $v_R=(1-\phi_R)u$ has finite semi-norm $[v_R]_{s,2}<\infty $, hence Proposition \ref{l:conv and Sobolev norm} gives $\rho_\eps*v_R\in C_0^\infty(\bG)\hookrightarrow \dom$ and
$\lim_{\eps\rightarrow 0}[v_R-\rho_\eps*v_R]_{s,2} =0.$ Thus $v_R\in\dom$ and $u$ is the limit of $\rho_\eps*v_R$ in the $[\, .\,]_{s,2}$ semi-norm.
\end{proof}

\subsection{The Rellich-Kondrachov compact embeddings}\label{ss:Relich}
We turn to a version of the compactness of the embedding of the Dirchlet space in the natural Lebesgue spaces on  a bounded domain in the group.  A related result in the case of Orlicz spaces  on Carnot groups was proven recently in  \cite[Theorem 2.4]{CMSV21} where, unlike the setting here, the spaces are all subspaces of $L^2(\bG)$ with
$\bG$ a Carnot group.

While we are not going to expand on Nikol'skii-Besov spaces, the following lemma uses crucially that we are estimating $R_hu-u$, $(R_hu)(g)=u(gh)$,  and the $\mathcal{D}^{s,2}(\bG)$ norm is invariant under left-translations since it involves left-translation invariant distance $d(g,h)=|h^{-1}g|$. If we were working with the right-invariant  $\mathcal{D}^{s,2}(\bG)$ space we could estimate  $L_hu-u$, where $(L_hu)(g)=u(hg)$.

\begin{lemma}\label{l:translation difference}
There exists a constant $C=C(Q,s)$ such that for $u\in L^{2}_{loc}(\bG) $ with $[u]_{s,2}<\infty$ and $(R_hu)(g)=u(gh)$ we have the inequality
\begin{equation}\label{e:translation difference}
 \norm {R_hu -u}_{L^2(\bG)}\leq C |h|^{s}[u]_{s,2}.
 \end{equation}
In particular, \eqref{e:translation difference} holds true for $u\in\dom$.
\end{lemma}

\begin{proof}
Let $r=|h|$ and $y\in B_r(g)=\{ y\in\bG\mid \, |y^{-1}g| < r\}$, cf. \eqref{e:right-invarinat balls}. By the triangle inequality  we have
\[
|u(gh)-u(g)|\leq |u(gh)-u(y)| + |u(y)-u(g)|.
\]
Averaging over the ball $B_r(g)$ we have
\[
|u(gh)-u(g)|\leq \fint_{B_r(g)}|u(gh)-u(y)|dy + \fint_{B_r(g)}|u(y)-u(g)|dy,
\]
hence using $(a+b)^2\leq 2(a^2+b^2)$ and Holder's inequality we have
\begin{equation}\label{e:triangle ineq for Lhu-u}
\frac 12 |u(gh)-u(g)|^2\leq \fint_{B_r(g)}|u(gh)-u(y)|^2dy + \fint_{B_r(g)}|u(y)-u(g)|^2 dy\equiv A(g) + B(g).
\end{equation}
Next, we estimate separately the integrals of $A$ and $B$. First, we compute
\begin{multline}
A(g)=\sigma_Q r^{-Q}\int_{|y^{-1}g|<r} |u(gh)-u(y)|^2dy = \sigma_Q r^{-Q}\int_{|hz|<r} |u(gh)-u(ghz)|^2dz\\
\leq  \sigma_Q r^{-Q}\int_{|z|<2r} |u(gh)-u(ghz)|^2dz
\end{multline}
by first letting $y=ghz$, hence $y^{-1}g=(hz)^{-1}$, and then noting that for $|y^{-1}g|<r$ we have from the triangle inequality
$
|z|=|h^{-1}hz|\leq |hz| +|h^{-1}|=|hz| +|h|< 2r.
$
Integrating over $\bG$, applying Fubini's theorem and using $|z|<2r$ we come to
\begin{multline}\label{e:A term}
\int_{\bG} A(g)\volg \leq \sigma_Q r^{-Q}\int_{\bG}\int_{|z|<2r} |u(gh)-u(ghz)|^2dz\volg\\
\leq \sigma_Q 2^{Q+2s} r^{2s}\int_{|z|<2r} \int_{\bG}\frac {|u(gh)-u(ghz)|^2}{|z|^{Q+2s}}\volg dz \leq  \sigma_Q 2^{Q+2s} r^{2s}\int_{\bG} \int_{\bG}\frac {|u(g)-u(gz)|^2}{|z|^{Q+2s}}\volg dz\\
\leq \sigma_Q 2^{Q+2s} r^{2s} [u]_{s,2},
\end{multline}
using the change of variable $g'=gh$ in the next to last inequality and, finally, recalling the definition of the $[\, . \,]_{s,2}$ semi-norm from \eqref{e:def frac norm}.

Turning to the estimate of the integral $B(g)$ defined in \eqref{e:triangle ineq for Lhu-u}, an integration over $\bG$ gives
\begin{multline}\label{e:B term}
\int_{\bG} B(g)\volg\leq \sigma_Q r^{-Q}\int_{\bG}\int_{|y^{-1}g|<r} |u(g)-u(y)|^2dy\volg \\
\leq \sigma_Q r^{2s}\int_{\bG}\int_{|y^{-1}g|<r} \frac {|u(g)-u(y)|^2}{r^{Q+2s}}dy\volg
\leq \sigma_Q r^{2s}\int_{\bG}\int_{|y^{-1}g|<r} \frac {|u(g)-u(y)|^2}{|y^{-1}g|^{Q+2s}}dy\volg\\
\leq \sigma_Q r^{2s}\int_{\bG}\int_{\bG} \frac {|u(g)-u(y)|^2}{|y^{-1}g|^{Q+2s}}dy\volg= \sigma_Q r^{2s}[u]_{s,2}.
\end{multline}
Finally, the claim of the Lemma follows from \eqref{e:triangle ineq for Lhu-u}, \eqref{e:A term} and \eqref{e:B term}.
\end{proof}

\begin{rmrk} A natural question is whether we can show with a constant independent of $s$ the inequality $$\norm {R_hu -u}_{L^2(\bG)}\leq C(1-s)^{1/2} |h|^{s}\norm{u}_{\mathcal{D}^{s,2}(\bG)}?$$ Using example  \ref{ex:non carnot}, the lack of the H\"ormander's bracket generating condition allows for a splitting of the variables, which leads to a negative answer as explained below.
\end{rmrk}

Let us consider the homogeneous group in Example \ref{ex:non carnot}. If the above inequality were true then taking the limit $s\rightarrow 1-$ will give with $\zeta=(\xi,\eta)$ the inequality
\[
\left ( \int_{\R^2} |u(x+\xi,y+\eta) - u(x,y)|^2{dxdy}\right)^{1/2} \leq  C_1\left( |\xi|+|\eta|^{1/2} \right) \left(\int_{\R^2}|\partial_x u (x,y)|^2 \, dxdy \right )^{1/2}
\]
For $u=\phi(x)\psi(y)$ this would imply
\begin{multline*}
\left ( \int_{\R^2} |\phi(x+\xi)\psi(y+\eta) - \phi(x)\psi(y)|^2{dxdy}\right)^{1/2}\\
 \leq  C_1\left( |\xi|+|\eta|^{1/2} \right)\left ( \int_{\R} |\partial_x\phi(x)|^2{dx}\right)^{1/2}\left ( \int_{\R} |\psi(y)|^2{dy}\right)^{1/2},
\end{multline*}
which is not true. For example, if we fix $\phi$ and take $\xi=0$ the above inequality reduces to
\[
\left ( \int_{\R} |\psi(y+\eta) - \psi(y)|^2{dy}\right)^{1/2}\\
 \leq  C_2 |\eta|^{1/2}\left ( \int_{\R} |\psi(y)|^2{dy}\right)^{1/2},
\]
which is false.

With the help of Lemma \ref{l:translation difference} and a Vitali covering argument we can prove the Rellich-Kondrachov type compact embedding Theorem \ref{t:frac Rellich embed}. We give a  short self-contained  proof, which can be considered as  a variation of the proof of Kolmogorov \cite{Kolm31}, see also \cite{Riesz33}, \cite[Lemma 11.1 and Theorem 11.1] {KV61} and \cite[Theorem 5]{H-OH10}.

\begin{thrm}\label{t:frac Rellich embed}
Let $\Omega$ be a bounded domain.  If $1\leq q< 2^*(s)$ then the inclusion $\iota:\mathcal{D}^{s,2}(\bG) \hookrightarrow L^q(\Omega)$ is compact.
\end{thrm}

\begin{proof}
First, we will prove the claim for $q=2$.  Thus, we will show that the restrictions to $\Omega$ of a bounded set in  $\mathcal{D}^{s,2}(\bG)$ is a totally bounded set in $L^2(\Omega)$. Let $\delta>0$. Consider the family of closed compact balls $$\mathcal{B}= \{\overline{B}_r(g)\mid r<\delta, g\in\Omega\}.$$
By the  proof of Vitali's covering theorem \cite[Theorem 1.6]{Hein01}
there exist finitely many disjoint balls  $\overline{B}_j=\overline{B}_{r_j}(g_j)\in \mathcal{B}$, $j=1, \dots, N$, that cover $\Omega$ up to a set of small measure, i.e., given $\eps>0$ we can assume that  $|\Omega\setminus \cup_{j=1}^N \overline{B}_j|<\eps$.
For  $u\in \mathcal{D}^{s,2}(\bG)$ let
\[
a_j=\fint_{\overline{B}_j}u(g) dg,  \quad    \text{and}\quad Pu(g)=\sum_{j=1}^N a_j\chi_j(g),
\]
where $\chi_j(g)$ is the characteristic function of the ball $\overline{B}_j$. Notice that for any $g$ the sum defining $Pu$ contains at most one non-zero term since the balls are pairwise disjoint. We will estimate the $L^2(\Omega)$ norm  of the difference $u-Pu$. With $\Omega_\eps=\Omega\setminus \cup_{j=1}^N \overline{B}_j$, $|\Omega_\eps|<\eps$, we split the integral into two integrals, which will be estimated separately
\begin{equation}\label{e:projection Vitali 0}
\norm {u-  Pu}^2_{L^2(\Omega)}= \int_{\Omega_\eps} |u(g)-Pu(g)|^2 \, dg \ + \ \sum_{i=1}^N \int_{\overline{B}_i}|u(g)- Pu(g)|^2 dg.
\end{equation}
The integral over $\Omega_\eps$ can be estimated using the  H\"older and the fractional Sobolev inequalities
\begin{multline}\label{e:projection Vitali 1}
\int_{\Omega_\eps} |u(g)-Pu(g)|^2 \, dg = \int_{\Omega_\eps} |u(g)|^2 \, dg \leq |\Omega_\eps|^{1 - \frac {2}{2^*(s)}} \norm{u}^2_{L^{2^*(s)}(\Omega_\eps)} \\
\leq  |\Omega_\eps|^{1 - \frac {2}{2^*(s)}} \norm{u}^2_{L^{2^*(s)}(\bG)} \leq C\eps^{\frac {2(Q-2s)s}{Q}} \norm{u}^2_{\mathcal{D}^{s,2}(\bG)}.
\end{multline}
Next, we estimate the last term in \eqref{e:projection Vitali 0} using that the balls $\overline{B}_j$ are pairwise disjoint, hence
\begin{equation}
\sum_{i=1}^N \int_{\overline{B}_i}|u(g)- Pu(g)|^2 dg = \sum_{i=1}^N \int_{\overline{B}_i}|u(g)- \sum_{j=1}^\infty a_j\chi_j(g)|^2 dg
= \sum_{i=1}^N \int_{\overline{B}_i}|u(g)-  a_i|^2 dg,
\end{equation}
where the last equality is obtained by noting that for $g\in \overline{B}_i$ we have $\chi_j(g)=0$ for $j\not=i$.
Therefore,
\begin{multline}\label{e:projection Vitali 2}
\sum_{i=1}^N \int_{\overline{B}_i}|u(g)- Pu(g)|^2 dg = \sum_{i=1}^N \int_{\overline{B}_i} \left\vert \fint_{\overline{B}_i}u(g)-  u(h)\, dh \right\vert^2 dg\\
 \leq \sum_{i=1}^N \int_{\overline{B}_i} \fint_{\overline{B}_i}\left\vert u(g)- u(h)\right\vert^2\, dh  dg,
\end{multline}
after applying H\"older's inequality to obtain the last  inequality. Next, we apply the change of variables $h=gh'$ in the inner integral and then use that for $g,h\in \overline{B}_j$ we have
\[
|h'|=|g^{-1}h|=|g^{-1}g_jg_j^{-1}h| \leq |g_j^{-1}g|+ |g_j^{-1}h|   \leq 2r_j\leq 2\delta,
\]
hence $h'$ belongs to the ball centered at the origin of radius $2r_j$, $h'\in \overline{B}_{2r_j}$. Therefore, we come to the estimate (denoting  $h'$ with $h$)
\begin{multline}\label{e:projection Vitali 3}
\sum_{i=1}^N \int_{\overline{B}_i}|u(g)- Pu(g)|^2 dg\leq \sum_{i=1}^N \frac {1}{|\overline{B}_{r_i}|}\int_{\overline{B}_i} \int_{\overline{B}_{2r_i}}\left\vert u(g)- u(gh)\right\vert^2\volh\volg\\
 \leq \sum_{i=1}^N \frac {1}{|\overline{B}_{r_i}|}\int_{\overline{B}_{2r_i}} \int_{\overline{B}_i} \left\vert u(g)- u(gh)\right\vert^2\volg\volh
\leq \sum_{i=1}^N \frac {1}{|\overline{B}_{r_i}|}\int_{\overline{B}_{2r_i}} \int_{\bG} \left\vert u(g)- u(gh)\right\vert^2\volg\volh \\
\leq
|| R_hu-u||^2_{L^2(\bG)}\sum_{i=1}^N \frac {|\overline{B}_{2r_i}|}{|\overline{B}_{r_i}|}
\leq C \delta^{2s}\norm{u}^2_{\mathcal{D}^{s,2}(\bG)},
\end{multline}
taking into account Lemma \ref{l:translation difference}, with $C$  a constant depending on $N$, $Q$ and $s$.

From \eqref{e:projection Vitali 1} and \eqref{e:projection Vitali 3}  it follows that for $\eps=\delta^{2^*(s)/2}$ we have
\begin{equation}\label{e:projection Vitali 4}
\norm {u-  Pu}^2_{L^2(\Omega)}\leq C \delta^{2s}\norm{u}^2_{\mathcal{D}^{s,2}(\bG)}.
\end{equation}
Furthermore, from H\"older's inequality the projection $P:L^2(\Omega)\rightarrow \R^N\cong \mathrm{span}\, \{ \chi_j \mid j=1,\dots, N\}$ is bounded.
Thus, given any $\epsilon >0$, by choosing a sufficiently small $\delta$ at the beginning of the proof  we can cover the image in $L^2(\Omega)$  of the given bounded set in $\mathcal{D}^{s,2}(\bG)$ with finitely many $\epsilon$-balls in $L^2(\Omega)$ centered at linear combinations of the functions $\chi_j(g)$, $j=1,\dots, N$, see  \cite[Lemma 1 and Theorem 5]{H-OH10}.

The case $q=2$ and H\"older's inequality imply the claim for $1\leq q\leq 2$ since $\Omega$ is bounded.  Finally, the case $2<q<2^*(s)$, hence $q=2(1-\theta)+\theta 2^*(s)$ for some $0<\theta<1$, follows from the interpolation inequality
\[
\norm {v}_{L^q(\Omega)} \leq \norm {v}^{1-\theta}_{L^2(\Omega)}\, \norm {v}^\theta_{L^{2^*(s)}(\Omega)},
\]
which, together with the already considered $L^2$-case, implies the sequential compactness of the embedding $\mathcal{D}^{s,2}(\bG)\hookrightarrow L^q(\Omega)$. This completes the proof of Theorem \ref{t:frac Rellich embed}.
\end{proof}

\begin{cor}\label{l:multiplication operator}
For a fixed $\phi\in C^\infty_0(\bG)$, let $M_\phi$ be the multiplication map defined on $\mathcal{D}^{s,2}(\bG) $ by $M_\phi(u)=\phi\, u$.  The map $M_\phi:\mathcal{D}^{s,2}(\bG)\rightarrow \mathcal{D}^{s,2}(\bG) $ is continuous. Furthermore, if $\Omega\Subset \bG$ and  $1\leq q< 2^*(s)$, then $\iota\circ M_\phi:\mathcal{D}^{s,2}(\bG)\hookrightarrow L^q(\Omega)$ is compact.
\end{cor}

\begin{proof}
 Suppose $\phi$ is supported in the ball $B_{R}$, $ \mathrm{supp}\, \phi\subset B_{R}$.
From the  triangle inequality we have for $x,\, y\in\bG$
\[
|\phi(x)u(x)-\phi(y)u(y)|\leq |\phi(x)||u(x)-u(y)| +|\phi(x)-\phi(y)||u(y)|,
\]
hence
 the $\mathcal{D}^{s,2}(\bG)$ norm of $\phi u$ can be estimated as follows
\begin{multline}\label{e:multiplication operator}
\frac 12\int_{\bG}\int_{\bG} \frac {|\phi(x)u(x)-\phi(y)u(y)|^2}{|y^{-1}x|^{Q+2s}}\, dxdy\\
\leq \norm {\phi}^2_{L^\infty(\bG)}\norm{u}^2_{\mathcal{D}^{s,2}(\bG)} +  \int_{\bG\setminus B_{2R}}\int_{\bG} \frac {|\phi(x)-\phi(y)|^2|u(y)|^2}{|y^{-1}x|^{Q+2s}}\, dxdy\ +  \ \int_{B_{2R}}\int_{\bG} \frac {|\phi(x)-\phi(y)|^2|u(y)|^2}{|y^{-1}x|^{Q+2s}}\, dxdy\\
=\norm {\phi}^2_{L^\infty(\bG)}\norm{u}^2_{\mathcal{D}^{s,2}(\bG)} + A +B.
\end{multline}
Since $\phi=0$ on the complement of the ball $B_R$ and $|y^{-1}x|\geq |y|/2$ for $y|\geq 2R$ and $|x|\leq R$ we have
\begin{multline*}
A=\int_{\bG\setminus B_{2R}}|u(y)|^2 \int_{\bG} \frac {|\phi(x)|^2}{|y^{-1}x|^{Q+2s}}\, dxdy\leq \norm {\phi}^2_{L^\infty(\bG)}\int_{\bG\setminus B_{2R}}|u(y)|^2 \int_{B_R} \frac {1}{|y^{-1}x|^{Q+2s}}\, dxdy\\
\leq 2^{Q+2s}\norm {\phi}^2_{L^\infty(\bG)} \int_{\bG\setminus B_{2R}}|u(y)|^2 \int_{\bG\setminus B_{2R}} \frac {1}{|y|^{Q+2s}}\, dxdy\\
\leq 2^{Q+2s}R^{Q}\norm {\phi}^2_{L^\infty(\bG)} \norm{u}^2_{L^{2^*(s)}(\bG\setminus B_{2R})}\left ( \int_{B_R} \frac {1}{|y|^{\frac {(Q+2s)Q}{2s}}}\, dy \right)^{2s/Q}\leq C\norm {\phi}^2_{L^\infty(\bG)} \norm{u}^2_{L^{2^*(s)}(\bG\setminus B_{2R})},
\end{multline*}
where $C=C(Q,s)$ after using H\"{o}lder's inequality noting that the H\"{o}lder conjugate to $2^*(s)/2$ is $(2^*(s)/2)'=Q/(2s)$ and the identity \eqref{e:polar coord}. Now, the fractional Sobolev type inequality \eqref{e:fracsobolev} gives
\begin{equation}\label{e:multiplication operator 2}
A\leq C\norm {\phi}^2_{L^\infty(\bG)} \norm{u}^2_{\mathcal{D}^{s,2}(\bG)}
\end{equation}
with a constant $C$ dependant on $s$ and $C$.

On the other hand, the integral $B$ in \eqref{e:multiplication operator} can be bounded in the same manner as we bounded   \eqref{e:small g case}, which gives
\begin{multline*}
B= \int_{B_{2R}}|u(y)|^2 \int_{\bG} \frac {|\phi(x)-\phi(y)|^2}{|y^{-1}x|^{Q+2s}}\, dxdy
= \int_{B_{2R}}|u(y)|^2 (D_s\phi)^2(y)\,dy
\\
\leq C(\norm {\phi}^2_{L^\infty(\bG)}+ \sum_{j=1}^n\norm {X_j\phi}^2_{L^\infty(\bG)})\int_{B_{2R}}|u(y)|^2 \,dy.
\end{multline*}
An application of H\"{o}lder's inequality gives
\begin{multline}\label{e:multiplication operator 3}
B\leq C\big(\norm {\phi}^2_{L^\infty(\bG)}+ \sum_{j=1}^n\norm {X_j\phi}^2_{L^\infty(\bG)}\big)\, R^{2s} \norm{u}^2_{L^{2^*(s)}(\bG\setminus B_{2R})} \\
\leq C\big(\norm {\phi}^2_{L^\infty(\bG)}+ \sum_{j=1}^n\norm {X_j\phi}^2_{L^\infty(\bG)}\big)\, R^{2s} \norm{u}^2_{\mathcal{D}^{s,2}(\bG)} ,
\end{multline}
with a constant $C$ dependant on $s$ and $C$ after using the fractional Sobolev inequality \eqref{e:fracsobolev}.

Inequalities \eqref{e:multiplication operator}, \eqref{e:multiplication operator 2} and \eqref{e:multiplication operator 3} imply the first part of the Lemma. The second claim follows from the first part and Theorem \ref{t:frac Rellich embed}.

\end{proof}

\subsection{Extremals of the fractional Sobolev embedding}\label{ss:extremals}
By Theorem \ref{t:frac homog Sobolev and Rellich embed} we have the continuous embedding $\mathcal{D}^{s,2}(\bG)\hookrightarrow L^{2^*(s)}(\bG)$. In particular, we can assume that the constant $S$ in \eqref{e:fracsobolev} is the best constant, i.e., it equals the norm of the embedding,

\begin{equation}\label{var1}
S^{-2}=I \overset{def}{=}\ \inf\ \left\{ \norm{u}^2_{\mathcal{D}^{s,2}(\bG)} \ \mid \ u\in C^{\infty}_0 (\bG),\ \underset{\bG}{\int} \lvert u (g) \rvert^{2^{*}(s)} dg =1 \right\}.
\end{equation}
The difficulty of finding an extremal is due to the translation and dilation invariance of the involved functionals.

The proof that the infimum in \eqref{var1} is obtained, similarly to the local case, with the help of  P. L. Lions' concentration compactness method \cite[Theorem 1.1]{L3}. In the present fractional setting, some modifications are necessary, and we briefly indicate them below. We note that all minimizers
of \eqref{var1} are of constant sign. Indeed, the $\mathcal{D}^{s,2}(\bG)$  norm of $|u|$ is not greater than the norm of $u\in
\mathcal{D}^{s,2}(\bG)$,   due to the elementary inequality
$$\left\vert\ | u(g)|-|u(h)|\ \right\vert \leq |u(g)-u(h)|.$$
Since the sign in the above inequality is strict if $u(g)u(h)<0$,
the $\mathcal{D}^{s,2}(\bG)$ norm is decreased by replacing a function with its absolute value, unless the function does not change sign. Furthermore, possibly after  rescaling, the (nonnegative) extremals achieving the infimum  in \eqref{var1} satisfy the
equation
\begin{equation}\label{e:fracYamabe1}
\frlap u=u^{2^*(s)-1}.
\end{equation}
The sharp asymptotic decay of  non-negative solutions of \eqref{e:fracYamabe1} was established in \cite{GLV22}.
As usual, for $u\in \mathcal{D}^{s,2}(\bG)$  the equation $\frlap u =F$ means
\begin{equation}\label{e:weak formulation}
\mathcal{D}(u,\phi)=\mathcal{D}_s(u,\phi)\overset{def}{=}
\int_{\bG}\int_{\bG} \frac
{(u(g)-u(h))(\phi(g)-\phi(h))}{| h^{-1}\, g|^{Q+2s}}\,
{dg}dh=\int_{\bG}F(g)\phi(g){dg},
\end{equation}
for any $\phi\in C^\infty_0(\bG)$.
For $\phi\in C^{\infty}_0 (B_R)$ from \eqref{e:grad square decay} we have that the carr\'e du champ function $D_s \phi=\sqrt{\Gamma_s(\phi)}$, recall \eqref{e:hor gradient}, satisfies $D_s \phi\in \mathcal{C}_\infty(\bG)$, the latter denoting the space of continuous functions on $\bG$ which vanish at infinity, i.e., the continuous functions $\phi$ on $\bG$ such that for every $\epsilon>0$ there is a compact $K\subset\bG$ such that $|\phi(g)|< \epsilon$, $g\in\bG\setminus K$.  Recall that the space $\mathcal{M}(\bG)$  of all bounded regular Borel measures on $\bG$ is the dual space of the Banach space $\mathcal{C}_\infty(\bG)$.

\subsection{The second concentration compactness lemma}
The ''second concentration compactness lemma'' of P.L. Lions, \cite[Lemma I.1]{L3} takes then the following form, where  $D_s u=\sqrt{\Gamma_s(u)}$ is the carr\'e du champ function, recall \eqref{e:hor gradient}.

\begin{lemma}\label{ccl2}
Suppose $ u_n \rightharpoonup u $ weakly  in $ \dom $,  while $ d\mu_n = (D_s u_n )^2 \vol \rightharpoonup
d\mu$ and  $d\nu_n=\lvert u_n\rvert^{2\text{*}(s)}\vol \rightharpoonup d\nu $ weak-$*$ in measure, where $ d\mu$ and $ d\nu $ are bounded, non-negative Borel measures on $ \bG. $ There exist points $ {g_j} \in \bG $ and  real numbers ${\nu_j}\ \geq\ 0$, ${\mu_j}\ \geq\ 0,$  at most countably many different from zero, such that, with $ I $ denoting the constant in \eqref{var1}, we have
\begin{equation}\label{e:ccl2}
\begin{aligned}
 d\nu & =\lvert u\rvert^{2^{*}(s)} \,dH+\sum_{j}{\nu_j} \delta_{g_j},\\
  d\mu & \geq (D_s u )^2\vol + \sum_{j}{\mu_j} \delta_{g_j},\\
  I{{\nu_j}}^{2/{2^{*}(s)}} & \leq {\mu_j}.
\end{aligned}
\end{equation}
In particular,

\begin{equation}
\sum {{\nu_j}}^{2/{2^{*}(s)}}<\infty.
\end{equation}
\end{lemma}

\begin{proof}
The  proof   can be achieved  as in \cite[Lemma I.1]{L3}, see also \cite[Lemma 1.4.5]{IV11}, taking into account the following comments. The needed point-wise convergence of the sequence $\{u_n\} $ follows from the Rellich type embedding  Theorem \ref{t:frac Rellich embed}.
With the notation of \cite[Lemma I.1, p. 161]{L3} where $v_n=u_n-u$,  letting
 \begin{equation}\label{e:meaures in c-c lema proof}
  d\lambda_n=(D_s v_n)^2dH \quad\text{and}\quad  d\omega_n\overset{def}{=}\bigl(\lvert u_n \rvert^{2^*(s)}-\lvert u \rvert^{2^*(s)}\bigr)dH=\lvert  v_n \rvert^{2^*(s)}dH+o(1)
 \end{equation}
 we have  $o(1)\rightharpoonup 0$ in measure by the Br\'ezis-Lieb lemma. We assume, by passing to a subsequence if necessary, that  $d\lambda_n\rightharpoonup d\lambda$ and  $d\omega_n\rightharpoonup d\omega=d\nu-\lvert u \rvert^{2^*(s)}dH$, both convergences weakly in the sense of  measures for some non-negative measures $d\lambda$  and $d\omega.$ The argument from \cite[(25)]{L3} leading to the reverse H\"older inequality \cite[(26)]{L3} needs to be modified slightly as follows.
Let
\[
F(g,h)=\frac {u\left(gh \right)\left(\phi\left(gh\right)-\phi\left(g\right)\right)}{| h|^{Q+2s}} \quad\text{and}\quad H(g,h)=\frac {\phi(g)\left(u(gh)-u(g)\right)}{| h|^{Q+2s}}.
\]
Note, that both $F$ and $H$ are functions in $L^2(\bG\times\bG)$, which can be seen arguing similarly to the first part of proposition \ref{l:frac lapl def and decay}.
The  definition of  the carr\'e du champ function $D_s \phi$, see \eqref{e:hor gradient}, and Minkowski's inequality imply the following inequality for the fractional gradient of product of two functions
\begin{multline}\label{e:fractional derivative of product norm inequality}
\norm {u\phi}_{\mathcal{D}^{s,2}(\bG)} = \norm{D_s(u\phi)}_ {L^2(\bG)}=\norm{F-H}_{L^2(\bG\times\bG)}\leq \norm{F}_{L^2(\bG\times\bG)}+\norm{H}_{L^2(\bG\times\bG)}\\
=\left(\int_{\bG}\int_{\bG} \frac {|u\left(gh \right)| ^2\left|\phi\left(gh\right)-\phi\left(g\right)\right|^2}{| h|^{Q+2s}}\,dgdh \right)^{1/2} +  \left( \int_{\bG}\int_{\bG} \frac {|\phi(g)|^2\left|u(gh)-u(g)\right|^2}{| h|^{Q+2s}}\,dhdg \right)^{1/2} \\
=  \left(\int_{\bG}\int_{\bG} \frac {|u\left(gh \right)| ^2\left|\phi\left(gh\right)-\phi\left(g\right)\right|^2}{| h|^{Q+2s}}\,dgdh \right)^{1/2} +\norm {\phi D_s u}_{L^2(\bG)}\\
= \norm {uD_s\phi}_{L^2(\bG)}+\norm {\phi D_s u}_{L^2(\bG)},
\end{multline}
noting that the change of variables $gh=h'$ in the first integral in the next to last line and Fubini's theorem give
\begin{multline}
 \int_{\bG}\int_{\bG} \frac {|u\left(gh \right)| ^2\left|\phi\left(gh\right)-\phi\left(g\right)\right|^2}{| h|^{Q+2s}}\,dhdg   =  \int_{\bG}\int_{\bG} \frac {|u\left(h \right)| ^2\left|\phi\left(h \right)-\phi\left(g\right)\right|^2}{| g^{-1}h|^{Q+2s}}\,dhdg\\
= \int_{\bG}\int_{\bG} \frac {|u\left(h \right)| ^2\left|\phi\left(h \right)-\phi\left(g\right)\right|^2}{| g^{-1}h|^{Q+2s}}\,dgdh = \int_{\bG}  |u \left(h \right)| ^2  \left| D_s \phi (h) \right|^2\,dh.
\end{multline}

 H\"older's inequality and the point-wise convergence imply that for $\phi\in C^{\infty}_0 (B_R)$ we have \begin{equation}\label{e:0 limit for RSI}
\int_{\bG} v_n^2(g) (D_s\phi)^2(g)\, dg \rightarrow 0.
\end{equation}
 From \eqref{e:fractional derivative of product norm inequality} and  \eqref{e:0 limit for RSI} we obtain
 \begin{multline}\label{e:limit for RSI}
 \int_{\bG} \lvert\phi\rvert^{2^*(s)}d\omega = \underset{n\rightarrow\infty}{\lim}\int_{\bG}\lvert\phi\rvert^{ 2^*(s)}d\omega_n=\underset{n\rightarrow\infty}{\lim}\int_{\bG}\lvert v_n\phi\rvert^{2^*(s)}dH\\
   \leq  I^{-{2^*(s)}/2}\underset{n\rightarrow\infty}{\lim} \Bigl(\int_{\bG} |D_s(v_n\phi)|^{2}dH\Bigr)^{{2^*(s)}/2}
   =I^{-{2^*(s)}/2}\underset{n\rightarrow\infty}{\lim} \Bigl(\int_{\bG}\lvert\phi\rvert^2  (D_sv_n)^{2}dH\Bigr)^{{2^*(s)}/2}\\
    = I^{-{2^*(s)}/2} \Bigl(\int_{\bG}\lvert \phi\rvert^2 d\lambda \Bigr)^{{2^*(s)}/2}.
 \end{multline}
The  limit \eqref{e:limit for RSI} together with the Sobolev type embedding lead to the reverse H\"older inequality
\begin{equation*}
 \left (\int_{\bG}\lvert\phi\rvert^{2^*(s)}d\omega\right)^{1/{2^*(s)}}\leq  I^{-1/2}\left ( \int_{\bG}\lvert\phi\rvert^2 d\lambda \right )^{1/{2}} .
 \end{equation*}
Hence, from \cite[Lemma 1.2]{L3}, we obtain
 \begin{equation*}
 d\nu=\lvert u \rvert^{2^*(s)}dH + \underset{j \in J}{\sum} {\nu_j}\delta_{g_j}\text{ and }  d\lambda\geq I \underset{j \in J}{\sum} {{\nu_j}}^{2/2^*(s)}.
 \end{equation*}
 From the weak convergence  $u_n\rightharpoonup u $ in $\dom$ it follows $d\mu_n- d\lambda_n=[D_s v_n]^2\, dH +o(1)$, which combined with the inequality above gives the desired estimate for $d\mu$.
\end{proof}

\subsection{Proof of the existence of extremals}
\begin{thrm}\label{t:existence of extremals}
The norm of the embedding $\mathcal{D}^{s,2}(\bG)\hookrightarrow L^{2^*(s)}(\bG)$ is achieved.
\end{thrm}

\begin{proof}
The proof, sketched below,  relies on the concentration compactness method \cite[Theorem 1.1]{L3}. Some modifications are needed due to the non-local setting of the considered problem. As before,  $D_s u=\sqrt{\Gamma_s(u)}$ is the carr\'e du champ function defined in \eqref{e:hor gradient}.

Let $u_n$ be a minimizing sequence of the variational problem \eqref{var1}.  From the weak-$*$ compactness of the unit ball, taking subsequences if necessary, we can assume $ d\mu_n = (D_s u_n )^2(g) dg \rightharpoonup
d\mu$ and  $d\nu_n=\lvert u_n(g)\rvert^{2\text{*}(s)}dg \rightharpoonup d\nu $ weak-$*$ in measure, where $ d\mu$ and $ d\nu $ are bounded, non-negative Borel measures on $ \bG$. In addition, by the Rellich type compactness, we can also assume $u_n\rightarrow u$ a.e. with respect to the fixed Haar measure. All of these properties will be preserved taking possibly rescaled or translated subsequences of the original sequence $\{u_n\}$.

The goal, as usual, is to see that after  suitable translations and dilations of the sequence $u_n$  using \eqref{e:scaling} we can obtain a subsequence, denoted by $v_n$, which converges in $L^{2^{*}(s)}(\bG)$ to a function $v\in L^{2^{*}(s)}(\bG)$ such that $\int_\bG \lvert v (g) \rvert^{2^{*}(s)} dg =1$.

For $u\in\mathcal{D}^{s,2}(\bG) $ let us define the concentration function
\begin{equation}\label{e:concentration function}
Q_u(r)=\sup_{g\in \bG}\int_{B_r(g)}\rho_u(g)\, dg, \qquad \rho_u(g)=|u(g)|^{2^*(s)}+  ({D_s u })^2(g).
\end{equation}

Using the dilations \eqref{e:scaling}, the continuity and range properties  of $Q_u(r)$, and the fact that $Q_{u_\lambda}(r)=Q_{u}(r/\lambda)$ we can replace the original minimizing sequence $u_n$ with a suitable  ''dilated'' and translated  sequence $v_n$ so that with $Q_n(r)=Q_{v_n}(r)$  we have
\begin{equation}\label{e:concentration function minimizer}
Q_n(1)={\int_{B_1(e)}} |v_n(g)|^{2^*(s)}\volg=1/2.
\end{equation}
By the choices made so far we also have, recall \eqref{var1},
\begin{equation}
 1+I=
 \lim_{n\rightarrow\infty}\int_{\bG}\rho_n(g)\, dg, \qquad \rho_n(g)\overset{def}{=}\rho_{v_n}(g)=|v_n(g)|^{2^*(s)}+  (D_sv_n)^2(g).
\end{equation}

In the next step of the proof, we apply the ''first concentration compactness lemma'' of P. L. Lions,  \cite[Lemma I.1]{L1} to the sequences of measures $\rho_n(g)dg$. The vanishing case is ruled-out by the normalization $Q_n(1)=1/2$.
Let $\rho_n(g)\, dg =d\nu_n +  d \mu_n$, where
\begin{equation}\label{e:prob measures}
d\nu_n=|v_n(g)|^{2^*(s)} \,  dg \qquad\text{and}\qquad d \mu_n= (D_sv_n)^2(g) \, dg.
\end{equation}
From the definitions we have $ v_n \rightharpoonup v $ weakly  in $ \mathcal{D}^{s,2}(\bG) $, while  $ d\mu_n = (D_s v_n)^2(g)\, dg \rightharpoonup
d\mu$ and  $d\nu_n=\lvert v_n\rvert^{2^*(s)}dg \rightharpoonup d\nu $ in the  weak-$*$ convergence in measure.
 The proof of \cite[Theorem 3.1]{MaMo19} carries over to our setting, hence  the dichotomy case cannot hold, as well. Thus,  the compactness case holds giving that the sequence $\rho_n(g)$ is tight, i.e., there exists a sequence of points $g_k\in\bG$, such that, for every $\eps>0$ there exists $R=R(\eps)$ such that for all sufficiently large $n$ we have
\begin{equation}\label{e:tightness 1}
\int_{ \bG\setminus B_{R}(g_k)} \rho_{n_k}(g_k)\, dg\leq \eps.
\end{equation}
Working with the subsubsequence $\rho_{n_k}$ we can assume
\begin{equation}
\underset{B_{R}(g_n)}{\int} d\nu_n \geq 1-\epsilon.
\end{equation}

If $\epsilon<1/2$ then $\underset{B_{R}(g_n)}{\int} d\nu_n >1/2$. Since by construction $\underset{B_{1}(e)}{\int} d\nu_n=1/2$ while $\underset{\bG}{\int} d\nu_n=1$ we see that $ {B_{1}(e)}$ and $ {B_{R}(g_n)}$ have a non-empty intersection, hence, using the triangle inequality,  $ {B_{R}(g_n)}\subset B_{2R+1}(e)$. This implies
\begin{equation}
\int_{B_{2R + 1}(e)} d\nu_n \geq \int_{B_{R}(g_n)} d\nu_n \geq 1-\epsilon,
\end{equation}
hence, after translating each of the $v_n$'s by $g_n$, the tightness inequality \eqref{e:tightness 1} holds with $g_k\equiv e$ for every $k$. By taking $\epsilon \rightarrow 0$ it follows that
\begin{equation}
\int_{\bG} d\nu=1.
\end{equation}

Next we apply Lemma \ref{ccl2}  to  the measures  $d\mu_n$ and the probability measures $d\nu_n$. For the sequence $\{d\mu_n\}$, we have $d\mu_n \rightharpoonup d\mu$ and $ \int_{\bG} d\mu_n \rightarrow I$, hence $\int_{\bG} d\mu \leq I$. We shall prove that all ${\nu_j}$'s in \eqref{e:ccl2} vanish, hence $ \int_{\bG}{\lvert v \rvert}^{2^*(s)}\, dg=1$. Let $ \alpha\overset{def}{=}\int_{\bG} {\lvert v\rvert}^{2^*(s)}\, dg \leq 1.$ Since $\int_{\bG} d\nu\ =\ 1$ it follows $ \sum {\nu_j}\ =\ 1-\alpha.$ From $\int_{\bG} d\mu\ \leq\ I$ we have $$\int_{\bG} (D_sv_n)^2(g) \, dg\ \leq \ I-\sum {\mu_j}.$$ From \eqref{e:ccl2} we have

\begin{equation*}
I\ =\ I_1\geq \underset{\bG}{\int} {(D_s v )}^2 dg+\sum {\mu_j}\geq I_{\alpha}+\sum I {{\nu_j}}^{2/2^*(s)} \ \geq\ I_{\alpha}  + \sum I_{{\nu_j}},
\end{equation*}
where
\begin{equation}\label{varlam}
I_\alpha\ \overset{def}{=}\ \inf\ \left\{ \underset{\bG}{\int} (D_sv_n)^2(g) \, dg\ : \ u\in C^{\infty}_o (\bG),\ \underset{\bG}{\int} \lvert u (g)\rvert^{2^*(s)} \volg = \alpha \right\}.
\end{equation}
Since $I_\alpha = {\alpha}^{2/2^*(s)} I_1$, it follows that $I_\alpha$ is strictly sub-additive, i.e.,
\begin{equation}\label{subadd}
I_1 < I_\alpha + I_{1-\alpha}, \text{ for every } 0<\alpha<1.
\end{equation}
The strict sub-additivity \eqref{subadd} of $I_\alpha$ implies that exactly one of the numbers $\alpha $ and ${\nu_j}$ is different from zero. We claim that $\alpha=1$. Indeed, suppose that there is a ${\nu_j}\ =\ 1$ and $d\nu\ =\ \delta_{g_j}$. From the normalization $ Q_n(1)=1/2$, hence
\begin{equation}\label{nopointmass}
1/2\ \geq\ \int_{B_1(g_j)} {\lvert v_n \rvert}^{2^*(s)}\, dg\ \rightarrow\ \int_{B_1(g_j)} d\nu \ =\ 1,
\end{equation}
which is a contradiction. Thus, we proved $\qn{v}=\alpha=1$ and $v_n \rightarrow v$ in $L^{2^*(s)}(\bG)$, which shows that $v$ is a solution of the variational problem \eqref{var1}.
\end{proof}

\section{Convolutions and approximations}\label{ss:convolutions}
We finish with a section containing a technical result of general interest concerning convolution and the fractional spaces in the considered homogeneous groups. We will consider as in \cite[p. 15]{FS82} the ''right convolution'' on $\bG$ defined by the formula
\begin{equation}\label{e:left convolution}
u*v (g)=\int_{\bG} u(gh^{-1}) v(h)dh= \int_{\bG} u(h) v(h^{-1}g)dh.
\end{equation}
The use of this convolution is motivated by the choice of the left-translation invariant distance, cf. \eqref{e:right inv dusatnce}.  The convolution is well defined for functions satisfying the conditions of the Young's inequality, \cite[Proposition (1.18) and (1.19)]{FS82}, if $1\leq p<\infty$, $1<q,\, r<\infty$, and $p^{-1}+q^{-1}=r^{-1}+1$ then we have
\[
\norm {u*v}_{L^{r,\infty}(\bG)} \leq C_{p,q} \norm{u}_{L^{p}(\bG)} \, \norm{v}_{L^{q,\infty}(\bG)}.
\]
Moreover, if $p>1$ then $\norm {u*v}_{L^{r}(\bG)} \leq C_{p,q} \norm{u}_{L^{p}(\bG)} \, \norm{v}_{L^{q,\infty}(\bG)}.$

The following lemma proves a general result, but the reader should keep in mind the cases of $u\in L^2(\bG)$  or $u\in L^{2^*(s)}(\bG)$ which are useful for controlled approximation of Sobolev functions via smooth functions.  An application is given in theorem \ref{t:trancation thrm}.

\begin{prop}\label{l:conv and Sobolev norm}
Let $\rho\in C^\infty_0(\bG)$ with $\int_{\bG} \rho dh =1$ and $\rho_\eps(g)=\eps^{-Q}\rho(\delta_{1/\eps} g)$. For
$u\in L^2_{loc}(\bG)$  and  $[u]_{s,2}<\infty$, denoted by $u_\eps=\rho_\eps * u$, we have the inequality
\begin{equation}\label{e:conv and Sobolev norm}
 [u_\eps]_{s,2}\leq [u]_{s,2}.
\end{equation}
Furthermore, we  have $[u_\eps - u]_{s,2}\to 0$ as $\eps\rightarrow 0$.
\end{prop}

\begin{proof}
 First, we recall that $\delta_\eps$ is a one-parameter group of group automorphisms, hence,
\begin{equation}\label{e:inverses of dilation}
\delta_{1/\eps}^{-1}h=\delta_{\eps}h\quad \text{and}\quad (\delta_\eps \, h)^{-1}=\delta_\eps \, h^{-1}.
\end{equation}
Therefore,  we have
\[
u_\eps\,  (x)= \eps^{-Q}\int_\bG  \rho(\delta_{1/\eps} z) u(z^{-1}x)dz = \int_\bG  \rho( z) u((\delta_\eps z)^{-1}x)dz
=\int_\bG  \rho( z) u((\delta_\eps z^{-1})x)dz.
\]
As in the Euclidean case \cite[(15.8)]{Po16} we continue the proof using Jensen's inequality.  Thus, with the help of Jensen's inequality, the invariance by left-translations of the distance and of  the Haar measure,  the norm of the convolution can be estimated as follows
\begin{multline*}
 {[u_\eps]^2_{s,2}}\leq \int_\bG\int_\bG \rho(z)\int_\bG \frac {|u((\delta_\eps z^{-1})x) - u((\delta_\eps z^{-1})y)|^2}{|y^{-1}x|^{Q+2s}}dxdy dz\\
=\int_\bG\int_\bG \rho(z)\int_\bG \frac {|u(x) - u(y)|^2}{|y^{-1}x|^{Q+2s}}dxdy dz=\int_\bG\int_\bG \frac {|u(x) - u(y)|^2}{|y^{-1}x|^{Q+2s}}dxdy \int_\bG \rho(z)  dz=[u]^2_{s,2}.
\end{multline*}

To prove the last part of the proposition, following \cite[Theorem 6.62]{Leoni23}, let $\triangle_h u$ be the function
$$\triangle_h u (g) = (R_hu-u)(g)=u(gh)-u(h).$$ From  lemma \ref{l:translation difference}  we have $\triangle_h u\in L^2(\bG)$.  From \eqref{e:def frac norm} we have
\begin{equation}\label{e:Rh and sp norm}
[{u}]^2_{{s,2}}
= \left(\int_{\bG}\int_{\bG}\frac {|\triangle_h u (g)|^2}{|h|^{Q+2s}}\, {dg}dh\right)^{1/2}= \left(\int_{\bG}\frac {\norm{\triangle_h u}_{L^2(\bG)}^2}{|h|^{Q+2s}}\, dh\right)^{1/2}<\infty. \end{equation}
The right-translation and mollification commute, i.e,  $R_hu_\eps=(R_hu)_\eps$, since
\[
    R_hu_\eps\,  (x)=u_\eps\,  (xh)= \eps^{-Q}\int_\bG  \rho(\delta_{1/\eps} z) u(z^{-1}x)dz = \eps^{-Q}\int_\bG  \rho(\delta_{1/\eps} z) (R_hu)(z^{-1}x)dz=(R_hu)_\eps\,  (x).
\]
Thus we have
\begin{equation}\label{e:[triangle,Rh]}
  (\triangle_h u)_\eps=\triangle_h u_\eps.
\end{equation}
From  Minkowski's inequality and the invariance of the measure by (left-)translations  it follows similarly to the proof of \eqref{e:conv and Sobolev norm} that mollification does not increase the $L^2$ norms, hence using \eqref{e:[triangle,Rh]} we obtain

\begin{equation}\label{e:Rh and mollification}
   \norm{\triangle_h u_\eps}_{L^2(\bG)}= \norm{(\triangle_h u)_\eps}_{L^2(\bG)}
    \leq\norm{\triangle_h u}_{L^2(\bG)}.
\end{equation}
 From \eqref{e:Rh and sp norm} we have
$$[u_\eps - u]_{s,2}=\left(\int_{\bG}\frac {\norm{\triangle_h (u_\eps - u)}_{L^2(\bG)}^2}{|h|^{Q+2s}}\, dh\right)^{1/2}. $$
By the properties of the mollification we have  $\norm{(\triangle_h u)_\eps- \triangle_h u}_{L^2(\bG)}\rightarrow 0$.
The claim follows from Lebesgue's dominated convergence theorem taking into account that $$\norm{\triangle_h (u_\eps - u)}_{L^2(\bG)}\leq 2 \norm{\triangle_h  u}_{L^2(\bG)}\in L^1({{|h|^{-Q-2s}}\, dh})$$ by \eqref{e:Rh and sp norm} since $[u]_{s,2}<\infty$.

\end{proof}

\end{document}